\documentclass[11pt]{article}

\RequirePackage[OT1]{fontenc}
\RequirePackage[utf8]{inputenc}
\RequirePackage{amsthm,amsmath,amsfonts}
\RequirePackage[numbers]{natbib}
\RequirePackage[colorlinks,citecolor=blue,urlcolor=blue]{hyperref}
\RequirePackage{authblk}

\newtheorem{Th}{Theorem}[section]

\newtheorem{Lmm}[Th]{Lemma} 
\newtheorem{Hy}[Th]{Hypothesis} 
\newtheorem{Defn}[Th]{Definition}

\title{Time Homogeneous Diffusion with Drift and Killing to Meet a Given Marginal }  
 
\author{John M. Noble \thanks{email address: {\tt noble@mimuw.edu.pl}\\accepted for publication: Stochastic Processes and their Applications  }}
\affil{Mathematical Statistics,\\ Institute of Applied Mathematics and Mechanics,\\  Faculty of Mathematics, Informatics and Mechanics, \\  University of Warsaw,\\ ul. Banacha 2,\\ 02-097 Warszawa, Poland}

\date{}
 
\begin{document}
 \maketitle

\begin{abstract}
In this article, it is proved that for any probability law $\mu$ over $\mathbb{R}$ and a drift field $b: \mathbb{R} \rightarrow \mathbb{R}$ and  killing field $k : \mathbb{R} \rightarrow \mathbb{R}_+$ which satisfy hypotheses stated in the article  and a given terminal time $t > 0$,  there exists a string $m$, an $\alpha \in (0,1]$, an initial condition $x_0 \in \mathbb{R}$ and a process $X$ with infinitesimal generator $\left(\frac{1}{2}\frac{\partial^2}{\partial m \partial x} + b \frac{\partial}{\partial m} - \frac{\partial K}{\partial m}\right)$ where $k = \frac{\partial K}{\partial x}$ such that for any Borel set $B \in {\cal B}(\mathbb{R})$, 
\[\mathbb{P}\left  ( X_t \in B | X_0 = x_0 \right ) = \alpha \mu (B).\]  

\noindent Firstly, it is shown the problem with drift and without killing can be accommodated, after a simple co-ordinate change, entirely by the proof in~\cite{N1}. The killing field presents additional problems and the proofs follow the lines  of~\cite{N1} with additional arguments.
\end{abstract}


 \section{Introduction} 
 \subsection{Results and Method of Proof}
Let  $\mu$ be a probability measure over $\mathbb{R}$, $b : \mathbb{R} \rightarrow \mathbb{R}$ and $k : \mathbb{R} \rightarrow \mathbb{R}_+$ given drift and killing functions. Set

\begin{equation}\label{eqtildb} \widetilde{b}(x) = \left\{ \begin{array}{ll} b(x) & x \in  \mbox{suppt}(\mu)  \\ 0 & x \not \in  \mbox{suppt}(\mu) \end{array}\right. , \qquad B(x) = \left\{ \begin{array}{ll} \int_{[0,x]} \widetilde{b}(y) dy & x \geq 0 \\ - \int_{[x,0)} \widetilde{b}(y) dy & x < 0  \end{array}\right. 
\end{equation}

\noindent where $\mbox{suppt}(\mu)$ denotes the support of the measure $\mu$.  Let 

\begin{equation}\label{eqtildk} \widehat{k}(x) = \left\{ \begin{array}{ll} k(x) & x \in  \mbox{suppt}(\mu)  \\ 0 & x \not \in  \mbox{suppt}(\mu)   \end{array}\right. , \qquad  K(x) = \left\{ \begin{array}{ll}  \int_{[0,x]} \widehat{k}(y) dy & x \geq 0 \\ - \int_{[x,0)} \widehat{k}(y) dy & x < 0. \end{array}\right.
\end{equation}

\begin{Hy}[Hypothesis on drift $b$, killing field $k$ and measure $\mu$]\label{hybmu} The target probability measure, drift and killing $(\mu, b, k)$ satisfy the following conditions. 
\begin{enumerate} \item $B$ from~\eqref{eqtildb} and $K$ from~\eqref{eqtildk} are absolutely continuous with respect to $\mu$. 
 \item \label{hycond2} Let $l_-(x) = \sup\{y \in \mbox{suppt}(\mu) \cap (-\infty, x)\}$ and let $l_+(x) = \inf\{y \in \mbox{suppt}(\mu) \cap (x,+\infty)\}$, then  

\begin{equation}\label{eqhypb2} \sup_{x \in \mathbb{R}}\lim_{h \downarrow 0} \int_{l_-(x) - h}^{l_+(x)+h} |\widetilde{b}(x)| dx < 1  
\end{equation}

\noindent where $\widetilde{b}$ is from~\eqref{eqtildb}.

\item \label{hycond3} Let $c: (0,1) \rightarrow \mathbb{R}_+$ denote the function defined by:

\begin{equation}\label{eqcgamm} c(x) =    \frac{\left ( \ln \frac{1}{x} \right ) - (1 - x)}{(1 - x)^2}.
\end{equation}

\noindent Let $\gamma$ satisfy:
\begin{equation}\label{eqgamdef} \gamma = \frac{1}{2} \left( 1 - \sup_{x \in \mathbb{R}} \lim_{h \downarrow 0} \int_{l_-(x) - h}^{l_+(x)+h} |\widetilde{b}(x)| dx \right ).
\end{equation} 

\noindent   Then  $(b,\mu)$  satisfies: 

\begin{equation}\label{eqhypb} 
  \int_{-\infty}^\infty \left( \int_{0 \wedge x}^{0 \vee x} e^{F (b,y)} dy   \right) \mu(dx)  < +\infty
 \end{equation}
 \noindent  where

\begin{equation}\label{eqefffF}  F (b,y) =    2\left ( \int_{0 \wedge y}^{0 \vee y} \left | \widetilde{b}(x)\right | dx + c(\gamma ) \sup_{\underline{t} :   (0\wedge y) = t_0 < \ldots < t_n = (0 \vee y)}\sum_{i=0}^{n-1} \left\{ \left(\int_{t_i}^{t_{i+1}} |\widetilde{b}(x)| dx\right)^2  \right\} \right ) \end{equation}

\noindent and $\widetilde{b}$ is defined by~\eqref{eqtildb}.  Here the maximum is taken over sequences of length $n$ for all $n \in \mathbb{N}$. 
\item \label{hypart4} $\lim_{x \rightarrow \pm +\infty} \frac{\partial K}{\partial \mu} (x) = 0$. 

Let $z_+ = \sup\{ x \in \mbox{suppt}(\mu)\}$ and $z_- = \inf\{x \in \mbox{suppt}(\mu)\}$. Then $\frac{\partial K}{\partial \mu}(x)$ is defined to be $0$ for $x > z_+$ and $x < z_-$. 
\end{enumerate}  
\end{Hy}

\noindent This article addresses the following problem: suppose that  $(\mu, b,  k)$ satisfy Hypothesis~\ref{hybmu}.  It is shown that there exists a string measure $m$, an $\alpha \in (0,1]$ and an $x_0 \in \mathbb{R}$ such that 

\begin{equation}\label{eqpr2}  \frac{1}{2}\frac{\partial^2}{ \partial m  \partial x} + \frac{\partial B}{\partial m} \nabla_m  - \frac{\partial K }{ \partial m }\end{equation}

\noindent where $\nabla_m$ is defined in Section~\ref{secdefigp} is the infinitesimal generator of a process $X$ satisfying 

\[ \mathbb{P}(X_t \in B | X_0 = x_0, X_t \not \in \{D\}) = \mu(B) \qquad \forall B \in {\cal B}(\mathbb{R})\]

\noindent where $D$ is a cemetery state, $X_t \in \{D\}$ denotes that the process has been killed by time $t$ and 

\[ \alpha := 1 - \mathbb{P}\left( X_t \in \{D\}\right) > 0.\]

\noindent If $t$ is replaced by an exponential time, $\alpha$, $x_0$ and $m$ are uniquely determined and an explicit construction is given. If $t$ is a deterministic time, only existence is given, although the method of proof may indicate how to provide approximations.

\paragraph{Remarks on Hypothesis~\ref{hybmu}} 

\begin{enumerate} 
\item \label{rq1} For $\gamma$ defined by~\eqref{eqgamdef}, it follows from~\eqref{eqhypb2} that $\gamma > 0$ (where the inequality is strict).
\item \label{rq2} For $x \in (0,1)$, the power series expansion of $1-x$ gives:
 \[ \log \frac{1}{x} = - \log (1 - (1 - x)) = \sum_{j=1}^\infty \frac{(1 - x)^j}{j}\]
 
 \noindent so that 
 
\[c(x) = \sum_{j=0}^\infty \frac{(1 - x)^{j}}{j+2}. \]

\noindent It follows that $\lim_{x \uparrow 1} c(x) = \frac{1}{2}$, $c(x)$ is decreasing in the range $x \in (0,1)$, $c(x) < +\infty$ for $x > 0$ and $\lim_{x \downarrow 0} c(x) = +\infty$. 

\item It is straightforward (and easier) to obtain the existence of a measure $m$  which gives an $\alpha > 0$ and a process with infinitesimal generator
\begin{equation}\label{eqprob2}
 \left(\frac{1}{2} \frac{\partial^2}{\partial m \partial x} + \frac{\partial B}{\partial m}\nabla_m  \right) - k
\end{equation}
\noindent for a given drift $b$ and killing $k$, which has distribution
\[ \mathbb{P}(X_\tau \in \{D\}) = 1-\alpha \qquad \mathbb{P}(X_\tau \in A) = \alpha \mu (A) \qquad \forall A \in {\cal B}(\mathbb{R}),\]

\noindent where $\tau$ is the terminal time, $\mu$ is the prescribed measure, $D$ denotes the cemetery state and $\{X_\tau \in \{D\}\}$ denotes that the process has been killed by time $\tau$. As with the case discussed in this article, with similar proofs, there is uniqueness and explicit construction when stopped at an independent geometric / exponential time. When finding a process with generator given by~\eqref{eqprob2}, the hypothesis on the killing field $k$ may be relaxed; Part~\ref{hypart4} of  Hypothesis~\ref{hybmu} is irrelevant for this problem, since it is only connected with ensuring that the limit of processes on atomised state spaces is not dead with probability 1 by the terminal time for a generator given by~\eqref{eqpr2}.  This issue resolves itself without this hypothesis for the generator given by~\eqref{eqprob2}.  
\end{enumerate}

\noindent The line of proof is as follows:

\begin{enumerate}
\item Discrete time and finite state space are considered; conditions under which a suitable Markov chain with a given distribution when stopped at an independent geometric time are established. The solution, when it  exists, is unique and the construction is explicit.  

\item This is then extended to establish conditions under which
there exists a Markov chain with a given distribution when stopped at an independent negative binomial time. This uses the fact that a negative binomial variable is the sum of independent identically distributed geometric variables and uses a fixed point theorem. For the problem of finding an infinitesimal generator of the form of~\eqref{eqpr2} or~\eqref{eqprob2}, substantial modifications of the arguments in~\cite{N1} are required when killing is introduced.   

\item Limits of negative binomial times   by reducing the time mesh are taken to obtain a time with Gamma distribution as in~\cite{N1}. Limits are then taken  to obtain a deterministic time. The arguments are along similar lines to those of~\cite{N1}, with some crucial modifications. 

\item Finally, arbitrary state space is considered. As in~\cite{N1}, the target
measure is approximated by a sequence of atomised measures. The drift is dealt with by a change of co-ordinates and the sequence of atomised measures in the transformed co-ordinates is considered. The killing is dealt with by considering the process without killing, together with the conditional distribution of the killing time. Both of these converge. The problem is to ensure that the diffusion coefficient does not tend to infinity and the probability that the process has been killed does not tend to $1$ as the limit is taken. The proof requires Hypothesis~\ref{hybmu} Part~\ref{hypart4}. 
\end{enumerate} 

\subsection{Background}
The problem of constructing a gap diffusion with a given law with compact support at an independent exponential time has been discussed fully by Cox, Hobson and
Ob\l{}\'{o}j in~\cite{CHO} (2011). The problem of constructing a martingale diffusion that has law $\mu$ at a fixed time $t$ has been solved by Jiang and Tao in~\cite{JT} (2001) under
certain smoothness assumptions. Recently, Forde in~\cite{F} (2011) extended the work of Cox, Hobson and Ob\l{}\'{o}j~\cite{CHO} to provide a process with prescribed joint law for the process at an independent exponential time $\tau$ and its supremum over the time interval $[0,\tau]$.   

For any prescribed measure $\mu$, the problem of finding a martingale diffusion with given marginal $\mu$ at a fixed time $t > 0$ was solved in~\cite{N1} (2013). Independently and simultaneously, Ekström, Hobson, Janson and Tysk~\cite{EHJT} (2013) found a different proof; in~\cite{EHJT}, the target distribution is again approximated by atomic measures, but general results from algebraic topology to conclude existence of a
limit. In~\cite{Mon} (1972), Monroe constructs a general symmetric stable process with a prescribed marginal at a fixed time, but does not require that the resulting process satisfies a martingale property.  

\subsection{Motivation}  The subject of strong Markov processes generated by Kre\u{i}n-Feller generalised second order differential operators  and, more specifically, the inverse problem of computing a function $a$ to give a solution 
$f$ to the parabolic equation

\[ \frac{\partial f}{\partial s} = a\left(\frac{1}{2} \frac{\partial^2}{\partial x^2} + b\frac{\partial}{\partial x} - k \right ) f \]

\noindent is of interest in its own right. Here $a$ is understood as $\frac{1}{m^\prime(x)}$ and the initial condition at $s = 0$ is a dirac mass $f(0,x) = \delta_{x_0}(x)$ at point $x_0 \in \mathbb{R}$; the end condition $f(t,x)$ for $s = t > 0$ is prescribed.   

The operator $\frac{\partial^2}{\partial m \partial x}$ and its spectral theory were introduced by Kre\u{i}n~\cite{Kr} (1952) and, for a more developed treatment, Kac and Kre\u{i}n~\cite{KacKr} (1958). A lucid account of the spectral theory is given by Dym and McKean~\cite{DMc} (1976). The operator, viewed as the generator of a strong Markov process, is discussed in Knight~\cite{Kn} (1981) where it is referred to as a {\em gap diffusion} and Kotani and Watanabe~\cite{KW} (1982) where it is referred to as a {\em generalised diffusion}. 

In recent years, interest in gap diffusion operators and their associated processes has been strongly renewed by applications to the field of modelling financial markets. The general motivating problem within finance is that of automating the pricing and risk management of derivative securities.  This is discussed by Carr and Nadtochiy in~\cite{Carr} (2014), where   the Local Variance Gamma model is developed to do this.   

The addition of drift and killing have importance when the prices of both the numéraire and the asset are modelled by stochastic processes. The covariation between the price of the numéraire and the price of the asset changes the drift of the discounted asset price process, hence the requirement to incorporate a drift $b$. The inclusion of a killing field extends the class of models available.

\paragraph{Acknowledgements} I thank Peter Carr for suggesting the problem of drift and killing and indicating the importance to financial applications. I also thank an anonymous referee whose thorough reading and careful comments led to substantial improvements.

\section{Definitions, Infinitesimal Generators and Processes}\label{secdefigp}

A definition of the operator ${\cal G} = \frac{\partial^2}{\partial m \partial x}$ used in   \eqref{eqpr2} may be found in Dym and McKean~\cite{DMc} or Kotani and Watanabe~\cite{KW}. The Kotani Watanabe definition is more useful in this setting, because it extends to strings defined over the whole real line. The domain of the operator, denoted ${\cal D}({\cal G})$ is the space of functions $f \in {\cal B}(\mathbb{R})$ such that there exists an $m$-measurable  function $g$ satisfying $\int_{x_1}^{x_2} g^2(x) m(dx) < +\infty$ for all $-\infty < x_1 < x_2 < +\infty$ such that 

\[ f(x) = f(x_0) + (x - x_0) f_-^\prime (x_0) + \int_{x_0}^x \int_{x_0 -}^y g(z) m(dz) dy \qquad \forall -\infty < x_0 < x < -\infty\]

\noindent where $f_-^\prime$ denotes the left derivative, $\int_a^b$ denotes integration over $(a,b]$ and $\int_{a-}^b$ denotes integration over the closed interval $[a,b]$. The quantity ${\cal G}f = \frac{\partial^2}{\partial m \partial x} f$ is defined as $g$. 

The operator $\nabla_m$ is defined as follows: let $z_- = \inf\{x \in \mathbb{R} | x \in \mbox{suppt}(m)\}$ and let $z_+ = \sup\{x \in \mathbb{R} | x \in \mbox{suppt}(m)\}$. For $x \in (z_-, z_+)$, define:
\[ \left\{\begin{array}{l} x^*_m (x) = \lim_{\epsilon \downarrow 0} \inf\{y > x+ \epsilon | y \in \mbox{suppt}(m)\} \\ x_{*m} (x)= \lim_{\epsilon \downarrow 0} \sup\{y < x-\epsilon | y \in \mbox{suppt}(m)\}. \end{array}\right. \]

\noindent The operator $\nabla_m$ is defined on functions $f \in {\cal D}({\cal G})$ as:

\begin{equation}\label{eqnabm}
 \nabla_m f(x) = \left\{ \begin{array}{ll} \lim_{h \rightarrow 0} \frac{f(x^*_m(x) + h) - f(x_{*m}(x) - h)}{x^*_m - x_{*m} + 2h} & x \in (z_-,z_+) \\ 0 & \mbox{other}. \end{array}\right. 
\end{equation}

\paragraph{Note 1} This definition of $\nabla_m$ is the definition associated with the drifts of the Markov chains under discussion. It boils down to Equation~\eqref{nabla} (given later) for a discrete state space and to $\frac{\partial}{\partial x}$ for $m(dx) = dx$.

\paragraph{Note 2} The definition of the domain of the operator is not discussed further in this article, since the method of proof does not require it, but it is reasonably straightforward to show that, for any  process $X$ obtained as the limit (in law) of processes with generators which converge to~\eqref{eqpr2} (as described in the article), if $f \in {\cal D}({\cal G})$, then for all $t > 0$, $F(t,.) \in {\cal D}({\cal G})$ where $F(t,x) := \mathbb{E}[f(X_t)|X_0 = x]$. It follows from the analysis given that $F$ thus defined satisfies:

\[ \frac{\partial}{\partial t} F = \frac{1}{2}\frac{\partial^2}{\partial m \partial x} F + \frac{\partial B}{\partial m} \nabla_m F - \frac{\partial K}{\partial m}F.\]

\noindent When $m$ has a well defined density $m^\prime > 0$, let $a = \frac{1}{m^\prime}$ then~\eqref{eqpr2} may be written as:

\begin{equation}\label{eqtwodens} a \left(\frac{1}{2}\frac{\partial^2}{\partial x^2} + b\frac{\partial}{\partial x} - k \right).\end{equation}

\noindent When a finite discrete state space ${\cal S} = \{i_1, \ldots, i_M\}$ is considered, the generator may be written as:

\[ a \left(\frac{1}{2}\Delta + b\nabla - k \right).\] 

\noindent where, for discrete state space, the definitions of the operators $\Delta$ and $\nabla$ are given in Definition~\ref{deflapder} below and, with abuse of notation, $b(i_j) : j= 2, \ldots, M-1$ here represents the sizes of the atoms of $B$ (Equation~\eqref{eqtildb}), with $b(i_1) = b(i_M) = 0$. 

\begin{Defn}[Laplacian and Derivative, Discrete state space]\label{deflapder}  Consider a state space
\[{\cal S} = \{i_1, \ldots, i_M\}, \qquad i_1 < \ldots < i_M.\]

\noindent For a function $f : {\cal S} \rightarrow \mathbb{R}$, the {\em Laplace operator} $\Delta$ is defined as: 

\begin{equation}\label{laplace} \left\{ \begin{array}{ll} \Delta f(i_1) = \Delta f(i_M) = 0 & \\
\Delta f(i_j) = \frac{2}{(i_{j+1} - i_j)(i_j - i_{j-1})}\left(\frac{i_j - i_{j-1}}{i_{j+1} - i_{j-1}} f(i_{j+1}) - f(i_j) + \frac{i_{j+1} - i_j}{i_{j+1} - i_{j-1}} f(i_{j-1})\right) & j = 2, \ldots, M-1 \end{array} \right.
\end{equation}

\noindent The {\em derivative operator} $\nabla$ is defined as: 

\begin{equation}\label{nabla} \left\{\begin{array}{ll} \nabla f(i_1) = \nabla f(i_M) = 0 & \\
\nabla f(i_j) = \frac{f(i_{j+1}) - f(i_{j-1})}{i_{j+1} - i_{j-1}} & 2 \leq j \leq M - 1 \end{array} \right. \end{equation} 
\end{Defn}

\paragraph{Remarks} \begin{enumerate} \item If the function $f$ is defined on an interval $(y_0,y_1)$, $f \in C^2((y_0,y_1))$ (twice differentiable with continuous second derivative) and a sequence ${\cal S}_n$ is considered, where ${\cal S}_n = \{i_{n,1}, \ldots, i_{n,M_n}\}$, $i_{n,j} < i_{n,j+1}$,  $i_{n,1} \downarrow y_0$, $i_{n,M_n} \uparrow y_1$ and $\lim_{n \rightarrow +\infty} \max_j (i_{n,j+1} - i_{n,j}) = 0$, with $\Delta_{{\cal S}_n}$ the operator defined on ${\cal S}_n$, then $\lim_{n \rightarrow +\infty} \Delta_{{\cal S}_n} f= \frac{d^2}{dx^2} f$. Note that the function $f$ has been defined on $C^2((y_0,y_1))$. The sense in which convergence is meant is: let $j_n(x) = \max\{ j : i_{n,j} \leq x\}$ then for all $x \in (y_0, y_1)$,

\begin{eqnarray*}\lefteqn{\lim_{n \rightarrow +\infty} \left | \frac{2}{(i_{n,j_n(x)+1} - i_{n,j_n(x)})(i_{n,j(x)} - i_{n,j(x)-1})}\left(\frac{i_{n,j(x)}    - i_{n,j(x)-1}}{i_{n,j(x)+1} - i_{n,j(x)-1}} f(i_{n,j(x)+1})\right. \right.} \\&& \left. \left. \hspace{20mm} - f(i_{n,j(x)}) + \frac{i_{n,j(x)+1} - i_{n,j(x)}}{i_{n,j(x)+1} - i_{n,j(x)-1}}f(i_{n,j(x)-1})\right) - \frac{d^2}{dx^2}f(x) \right | = 0.
\end{eqnarray*}

\item If the function $f$ is defined on the whole interval $(y_0,y_1)$ and $f \in C^1(\mathbb{R})$ (differentiable, continuous first derivative) and a sequence ${\cal S}_n$ is considered where  ${\cal S}_n = \{i_{1,n}, \ldots, i_{M,n}\}$, $i_{j,n} < i_{j+1,n}$  $i_{1,n} \downarrow y_0$ and $i_{M,n} \uparrow y_1$ and \[ \lim_{n \rightarrow +\infty} \max_j (i_{j+1,n} - i_{j,n}) = 0,\] with $\nabla_n$ the operator defined on ${\cal S}_n$, then $\lim_{n \rightarrow +\infty} \nabla_n f = \frac{d}{dx} f$ in the sense that for $f \in C^1((y_0,y_1))$ (differentiable with continuous derivative, for all $x \in (y_0,y_1)$
\[ \lim_{n \rightarrow +\infty} \left | \frac{f(i_{j_n(x)+1}) - f(i_{j_n(x)-1})}{i_{j_n(x)+1} - i_{j_n(x)-1}} - \frac{d}{dx}f(x) \right | = 0\]
if $\lim_{n \rightarrow +\infty} \max_j(i_{n,j+1} - i_{n,j}) = 0$.
\end{enumerate}

\paragraph{Notation} For finite discrete state space ${\cal S} = \{i_1, \ldots, i_M\}$,   $b(i_j)$ and $k(i_j)$ are used to denote the sizes of the {\em atoms} of $B$ and $K$ respectively from Equations~\eqref{eqtildb} and~\eqref{eqtildk}. This is a minor abuse of notation, since in~\eqref{eqtildb} and~\eqref{eqtildk}, $b$ and $k$ are used to denote the {\em derivatives} of $B$ and $K$ on $\mbox{suppt}(\mu)$. This notation will be used throughout when dealing with the problem on discrete state space. 

Furthermore, for finite discrete state space ${\cal S} = \{i_1, \ldots i_M\}$, let $b: {\cal S}\backslash\{i_1,i_M\} \rightarrow \mathbb{R}$ and $k :{\cal S}\backslash\{i_1,i_M\} \rightarrow \mathbb{R}_+$ denote the drift and killing respectively. The following notation will be used: let $\underline{b} = (b_2, \ldots, b_{M-1})$ where (with slight abuse of notation) $b_j = b(i_j)$ and $\underline{k} = (k_2, \ldots, k_{M-1})$ where (same notation) $k_j = k(i_j)$.   The notation $\widetilde{k}_j$ will be used to denote the following:

\begin{equation}\label{eqtildek} \left\{\begin{array}{ll} \widetilde{k}_j = (i_{j+1} - i_j)(i_j - i_{j-1}) k_j & j \in (2, \ldots, M-1)  \\ \widetilde{k}_1 = \widetilde{k}_M = 0 \end{array}\right. \end{equation}

\noindent and $\underline{\widetilde{k}} = (\widetilde{k}_1,  \ldots, \widetilde{k}_{M})$. \qed \vspace{5mm}

\noindent For all results with finite state space, the following hypothesis will be required:

\begin{Hy} \label{hyb} For a discrete, finite state space ${\cal S} = \{i_1, \ldots i_M\}$ where $i_1 < \ldots < i_M$, the vector $\underline{b}$ satisfies the condition:

\begin{equation}\label{eqbhyp} -\frac{1}{i_{j+1} - i_j}  < b_j <  \frac{1}{i_j - i_{j-1}} \qquad j = 2, \ldots, M-1.\end{equation}
\end{Hy}

\noindent Set

\begin{equation}\label{eqdefq} \left\{ \begin{array}{ll}
q_{j,j+1} = \frac{i_j - i_{j-1}}{i_{j+1} - i_{j-1}}(1 +  (i_{j+1} - i_j) b_j) & j = 2, \ldots, M-1 \\
q_{j,j-1} = \frac{i_{j+1} - i_j}{i_{j+1} - i_{j-1}}(1 - (i_j - i_{j-1}) b_j) & j = 2, \ldots, M-1
\end{array} \right. 
\end{equation}
           
\noindent Condition~\eqref{eqbhyp} is necessary and sufficient to ensure that $q_{j,j+1}$ and $q_{j,j-1}$ are non negative for each $j$. With these definitions of $\underline{b}$, $\underline{k}$ and $\underline{\widetilde{k}}$, the following definitions are made for the transitions (in discrete time) and the intensities (in continuous time) of the Markov processes that are of interest.

\begin{Defn}[Transition Matrix] \label{deftrmat} Let $\underline{\lambda} = (\lambda_2, \ldots, \lambda_{M-1}) \in \mathbb{R}_+^{M-2}$. For $h < \frac{1}{\max_{j \in \{2,\ldots, M-1\}} \lambda_j (1 + \widetilde{k}_j)}$, let $\widetilde{P}^{(h)}(\underline{k}, \underline{\lambda})$ be  the $M+1 \times M+1$ matrix defined by:

\begin{equation}\label{eqptild} \left\{ \begin{array}{ll} \widetilde{P}^{(h)}_{j,M+1} (\underline{k}, \underline{\lambda}) = h \lambda_j \widetilde{k}_j & j=2,\ldots, M-1 \\
\widetilde{P}^{(h)}_{1,M+1} (\underline{k}, \underline{\lambda})= \widetilde{P}^{(h)}_{M,M+1}(\underline{k}, \underline{\lambda}) = 0 & \\ \widetilde{P}^{(h)}_{M+1,M+1}(\underline{k}, \underline{\lambda}) = 1 & \\ \widetilde{P}^{(h)}_{M+1,j}(\underline{k}, \underline{\lambda}) = 0 & j = 1, \ldots, M \\  \widetilde{P}^{(h)}_{jj}(\underline{k}, \underline{\lambda}) = 1 - \lambda_j (1  +   \widetilde{k}_j)h & j = 2, \ldots, M-1 \\ \widetilde{P}^{(h)}_{11}(\underline{k},\underline{\lambda}) = \widetilde{P}^{(h)}_{MM}(\underline{k},\underline{\lambda}) = 1 & \\
 \widetilde{P}^{(h)}_{12} (\underline{k}, \underline{\lambda})= \widetilde{P}^{(h)}_{M,M-1} (\underline{k}, \underline{\lambda}) = 0 & \\ 
\widetilde{P}^{(h)}_{j,j+1}(\underline{k}, \underline{\lambda}) = h \lambda_j q_{j,j+1} & j = 2, \ldots, M-1\\ 
\widetilde{P}^{(h)}_{j,j-1}(\underline{k}, \underline{\lambda}) = h \lambda_j q_{j,j-1} & j = 2, \ldots, M-1  \\ \widetilde{P}_{jk}^{(h)} (\underline{k}, \underline{\lambda}) = 0 & |j-k| \geq 2, \quad (j,k) \in \{1, \ldots, M\}^2  \end{array} \right. \end{equation}
Let $P^{(h)}(\underline{k}, \underline{\lambda})$ denote the $M \times M$ matrix defined by $P^{(h)}_{ij}(\underline{k}, \underline{\lambda}) = \widetilde{P}^{(h)}_{ij}(\underline{k}, \underline{\lambda})$ for $(i,j) \in \{1, \ldots, M\}^2$.
\end{Defn}

\begin{Defn}[Intensity Matrix]\label{defintmat} Let 

\begin{equation}\label{eqthedef} \Theta(\underline{k},\underline{\lambda}) = \lim_{h \rightarrow 0}\frac{1}{h} \left(\widetilde{P}^{(h)}(\underline{k},\underline{\lambda}) - I \right)
\end{equation}
\end{Defn} 

\noindent It is straightforward to see that the matrix $\Theta (\underline{k}, \underline{\lambda})$ satisfies:

\begin{equation}\label{eqthetmat}
\left\{ \begin{array}{ll} \Theta_{j,M+1}(\underline{k}, \underline{\lambda}) = \lambda_j \widetilde{k}_j & j = 2, \ldots, M - 1\\
\Theta_{1,M+1} (\underline{k}, \underline{\lambda})= \Theta_{M,M+1}(\underline{k}, \underline{\lambda}) = 0 & \\ 
\Theta_{M+1,j}(\underline{k}, \underline{\lambda})= 0 & j = 1, \ldots, M+1\\
\Theta_{jj}(\underline{k}, \underline{\lambda}) = - \lambda_j ( 1 + \widetilde{k}_j) & j = 2, \ldots, M-1 \\ 
\Theta_{11}(\underline{k}, \underline{\lambda}) = \Theta_{MM}(\underline{k},\underline{\lambda}) = 0 & \\ 
\Theta_{12} (\underline{k}, \underline{\lambda}) = \Theta_{M,M-1}(\underline{k}, \underline{\lambda}) = 0 & \\
\Theta_{j,j+1}(\underline{k}, \underline{\lambda}) = \lambda_j q_{j,j+1} & j = 2,\ldots, M-1 \\
\Theta_{j,j-1}(\underline{k}, \underline{\lambda}) = \lambda_j q_{j,j-1} & j = 2, \ldots, M-1\\
\Theta_{j,k} (\underline{k}, \underline{\lambda}) = 0 & |j - k| \geq 2, \qquad (j,k) \in \{1, \ldots, M\}
\end{array}\right. 
\end{equation}

\noindent which is the {\em intensity matrix} of a Continuous Time Markov Chain on state space $\{1, \ldots, M+1\}$.\vspace{5mm}  

\noindent {\bf Note} The dependence on $\underline{k}$ and $\underline{\lambda}$ for $\widetilde{P}^{(h)}(\underline{k}, \underline{\lambda})$, $P^{(h)}(\underline{k}, \underline{\lambda})$ and $\Theta(\underline{k}, \underline{\lambda})$  will be suppressed; these will be written as $\widetilde{P}^{(h)}$, $P^{(h)}$ and $\Theta$ respectively.  \vspace{5mm}

\noindent For $h < \frac{1}{\max_{j \in \{2, \ldots, M-1\}}  \lambda_j ( 1 + \widetilde{k}_j)}$, $\widetilde{P}^{(h)}$  is the one-step transition matrix for a time homogeneous Markov process $X^{(h)}$,  with time step length $h$, satisfying

\[ \mathbb{P}(X^{(h)}_{h(t+1)} = i_k | X^{(h)}_{ht} = i_j) = \widetilde{P}^{(h)}_{jk}.\]

\noindent As discussed in~\cite{N1}, as $h \rightarrow 0$, the process $X^{(h)} \rightarrow X$ (convergence in the sense of finite dimensional marginals) to a continuous time Markov chain with intensity matrix $\Theta = \lim_{h \rightarrow 0}\frac{1}{h}(\widetilde{P}^{(h)} - I)$ from Definition~\ref{defintmat}.

\begin{Lmm} \label{lmmfininfgen2} Let ${\cal S} = \{i_1, \ldots, i_M\}$ and let $X$ be a continuous time Markov process on $S \cup \{D\}$ with transition intensity matrix from Definition~\ref{defintmat} Equation~\eqref{eqthedef} in the sense that

\[ \left\{ \begin{array}{ll} \lim_{h \rightarrow 0}\frac{1}{h} \mathbb{P}\left(Y_h = i_k | Y_0 = i_j\right) = \Theta_{jk} & (j,k) \in \{1, \ldots, M+1\}^2 \qquad j \neq k,\\ \lim_{h \rightarrow 0}\frac{1}{h}\left(\mathbb{P}\left(Y_h = i_j | Y_0 = i_j\right) - 1 \right ) = \Theta_{jj} & j \in \{1, \ldots, M+1\} \end{array}\right. \]

\noindent Let 
\begin{equation}\label{eqaj}\left\{\begin{array}{l}  a_j =  \lambda_j (i_{j+1} - i_j)(i_j - i_{j-1}) \qquad  j = 2, \ldots, M-1\\ a_1 = a_M = 0\end{array}\right. 
\end{equation} 

\noindent and (where the notation is clear) let $a : {\cal S} \rightarrow \mathbb{R}_+$ denote the function defined by $a(i_1) = a(i_M) = 0$, $a(i_j) = a_j$ for $j = 2, \ldots, M-1$. Then $Y$ has infinitesimal generator
\[ a \left(\frac{1}{2}\Delta + b \nabla - k\right).\]
\end{Lmm}

\paragraph{Proof} Recall the definition of $\widetilde{\underline{k}}$ (Equation \eqref{eqtildek}). Let $f$ be a function defined on $\{i_1, \ldots, i_M\}$ and let $F(t,i_j) = \mathbb{E}_{i_j}\left [ f(Y_t)\right ]$. Then, for $j = 2, \ldots, M-1$,

\begin{eqnarray*} \frac{\partial}{\partial t} F(t, i_j) &=& \lim_{h \rightarrow 0} \frac{F(t+h, i_j) - F(t, i_j)}{h} = \lim_{h \rightarrow 0} \frac{1}{h} \mathbb{E}_{i_j} \left [f(Y_{t+h}) - f(Y_t) \right] \\ &=& \lim_{h \rightarrow 0} \frac{1}{h}\left(\mathbb{E}_{i_j}[F(t,Y_h)] - F(t,i_j)\right)\\
&=& \lambda_j\left( q_{j,j+1} F(t,i_{j+1}) - F(t,i_j) + q_{j,j-1} F(t,i_{j-1}) - \widetilde{k}_j F(t,i_j)\right) \\
&=&  \lambda_j \left(\left(\frac{i_j - i_{j-1}}{i_{j+1} - i_{j-1}}F(t,i_{j+1}) - F(t,i_j) + \frac{i_{j+1} - i_j}{i_{j+1} - i_{j-1}} F(t,i_{j-1})\right)\right. \\&& + \left. \left( (i_{j+1} - i_j)(i_j - i_{j-1}) \right)b_j \left (\frac{F(t,i_{j+1}) - F(t,i_{j-1})}{i_{j+1} - i_{j-1}} \right )   - \widetilde{k}_j F(t,i_j)\right )\\
&=& a_j \left(\frac{1}{2}\Delta + b_j \nabla - k_j \right ) F(t,i_j) 
\end{eqnarray*}

\noindent For $j \in \{1, M\}$,

\[ \frac{\partial}{\partial t} F(t,i_j) = 0\]

\noindent as required. \qed

\section{Coordinate change to deal with the drift}\label{subcoc} The addition of the  drift $b$ can be dealt with through a simple change of co-ordinates, described here. The aim is to find a mapping of the process from space ${\cal S}$ (the state space of the process) to a space ${\cal R}$ such that the transformed process is drift free. For finite state space, ${\cal S} = \{i_1, \ldots, i_M\}$, $i_1 < \ldots < i_M$ the aim is to find a map $\kappa$ where (with abuse of notation) $\kappa_j = \kappa(i_j)$ for $j = 1, \ldots, M$ where $\kappa_j < \kappa_{j+1}$ for $j = 1, \ldots, M-1$ such that $q_{j,j-1}$ and $q_{j,j+1}$, for $j = 2, \ldots, M-1$, defined by~\eqref{eqdefq}, satisfy:

\begin{equation}\label{eqkap1} q_{j,j+1} = \frac{\kappa_j - \kappa_{j-1}}{\kappa_{j+1} - \kappa_{j-1}} \qquad q_{j,j-1} = \frac{\kappa_{j+1} - \kappa_j}{\kappa_{j+1} - \kappa_{j-1}}.
\end{equation}

\noindent Let 

\[ \delta_j = i_j - i_{j-1}, \qquad j = 2, \ldots, M.\]

\noindent Directly from~\eqref{eqdefq} and~\eqref{eqkap1}, it follows that $\kappa_{j} - \kappa_{j-1} = \epsilon_j$ for $j = 2, \ldots, M$ where  $\epsilon_2, \ldots, \epsilon_M$ satisfy

\[ \frac{\epsilon_j}{\epsilon_{j+1} + \epsilon_j} = \frac{\delta_j}{\delta_j + \delta_{j+1}} + \frac{\delta_j \delta_{j+1}}{\delta_j + \delta_{j+1}} b_j.\]

\noindent It follows that $\underline{\epsilon} = (\epsilon_2, \ldots, \epsilon_M)$ satisfies

\begin{equation}\label{eqeps} \frac{\epsilon_{j+1}}{\epsilon_j} = \frac{\delta_{j+1}}{\delta_j}\frac{1 - \delta_j b_j}{1 + \delta_{j+1}b_j} \qquad j = 2, \ldots, M-1.
\end{equation}

\noindent Clearly,~\eqref{eqeps} does not determine $\underline{\kappa}$ uniquely; two additional conditions have to be specified, which represent centring and scaling.  The following choice is made:  let 

\begin{equation}\label{eqemminemplus}
e_- = \inf \left \{j : \sum_{i=1}^j p_i \geq \alpha \right \}, \qquad e_+ = \sup \left \{j : \sum_{i= j}^M p_i \geq \alpha \right \}.
\end{equation}

\noindent where $0 < \alpha < 0.5$ is a number chosen such that   $e_-  < e_+$ (the inequality is strict). This is possible if ${\cal S}$ has $3$ or more distinct states. Let $K = (i_{e_+} - i_{e_-}) \vee 1$. Then $\underline{\epsilon} = (\epsilon_2, \ldots, \epsilon_M)$ defined by

\begin{equation}\label{eqeps2} \epsilon_{j+1} = K 
\frac{\delta_{j+1} \prod_{k=1}^j \left(\frac{1 - \delta_k b_k }{1 + \delta_{k+1} b_k}\right)}
{\sum_{a= e_-}^{e_+-1} \delta_{a+1} \prod_{k=1}^a \left(\frac{1 - \delta_k b_k}{1 + \delta_{k+1} b_k}\right ) }.
\end{equation}

\noindent satisfies~\eqref{eqeps} and $\kappa_{e_+} - \kappa_{e-} = K$. \vspace{5mm}

\noindent Define $\underline{\kappa}$ in the following way: let $e$  satisfy $\sum_{i=1}^e p_i \geq \frac{1}{2}$ and $\sum_{i=e}^M p_i \geq \frac{1}{2}$ and set $\kappa_e = 0$. For $j \neq e$, set

\begin{equation}\label{eqkappa} 
\left\{ \begin{array}{ll} \kappa_{e+j + 1} = \kappa_{e + j} + \epsilon_{e+j+1}  & j = 0, \ldots, M-e-1 \\ \kappa_{e - j -1} = \kappa_{e - j} - \epsilon_{e - j - 1} & j = 0, \ldots, e - 2 
\end{array}\right. 
\end{equation}

\noindent where $\underline{\epsilon}$ is defined by~\eqref{eqeps2}. It is easy to see that $\underline{\kappa}$, thus defined, satisfies~\eqref{eqkap1}.

This choice ensures that $\underline{\kappa}$ is centred (so that in the transformed coordinates the process is a martingale with mean zero)  and, furthermore, that when the discussion is extended in Section~\ref{seccoord} to the case of arbitrary measure $\mu$ on $\mathbb{R}$ with an appropriate sequence of atomised measures $\mu^{(N)}$, the processes with state spaces $\underline{\kappa}^{(N)}$  have suitable convergence properties.

\section{A Function to Accommodate the Killing Field} The method of proof adopted in this article is to try and rephrase the problem, as much as possible, in the language of~\cite{N1} and to use as much of the technique from~\cite{N1} as possible. The previous section introduced a co-ordinate change to deal with the drift $b$; under the co-ordinate change, the problem with drift, but without killing, reduces to that of~\cite{N1}. The introduction of killing presents other problems: firstly, even without drift, the initial condition is no longer as clear as it was in~\cite{N1} when the killing field is non trivial. It cannot be taken as simply the expectation of the target distribution, since the process in the time interval $[0,t]$, conditioned on being alive at time $t$, is no longer a martingale. Secondly,  the process is killed at rate $\lambda_j \widetilde{k}_j$ on site $i_j$, where $\underline{\lambda}$ is the holding intensity vector which is to be computed. This feeds into the equation required to obtain the intensities and there is no longer an explicit expression, even for the process stopped at an independent exponential time, like the formula $\lambda_j = \frac{1}{p_j}{\cal F}_j(\underline{p})$ that was available in~\cite{N1}. The function ${\cal G}$ described in this section plays the role of ${\cal F}$ in~\cite{N1}.

Let $\underline{p} = (p_1, \ldots, p_M)$ satisfy $\min_j p_j > 0$ and $\sum_{j=1}^M p_j = 1$. Let ${\cal S} = \{i_1, \ldots, i_M\}$, $i_1 < \ldots < i_M$ be a finite state space with $M$ elements; $\underline{i} = (i_1, \ldots, i_M)$ will be used to denote the elements of the space. Let $\underline{b} = (b_2, \ldots, b_{M-1})$ satisfy Hypothesis~\ref{hyb} and let $\underline{\kappa} = (\kappa_1, \ldots, \kappa_M)$ denote the coordinate change of $\underline{i}$ defined by Section~\ref{subcoc}. Let $\underline{k} = (k_2, \ldots, k_{M-1}) \in  \mathbb{R}_+^{M-2}$ denote the killing field and let $\widetilde{\underline{k}}$ be defined by Equation~\eqref{eqtildek}.  Let ${\cal G}_j(t,\underline{p}) : j \in \{ 1,\ldots, M\}$ be defined as follows: 

\begin{equation}\label{eqgeee} \left\{ \begin{array}{l} {\cal G}_1 = {\cal G}_M = 0 \\ {\cal G}_j  = \frac{1}{tp_j} \frac{(\kappa_{j+1} - \kappa_{j-1})}{(\kappa_{j+1} - \kappa_j)(\kappa_j - \kappa_{j-1})}\times \left\{  \begin{array}{ll} \sum_{a = 1}^{j-1} (\kappa_j - \kappa_a)p_a (1 + t{\cal G}_a \widetilde{k}_a) & 2 \leq j \leq l-1 \\ \sum_{a = j+1}^M (\kappa_a - \kappa_j) p_a (1 + t {\cal G}_a \widetilde{k}_a) & l \leq j \leq M-1 \end{array}\right.  \\ l : \kappa_{l-1} < \frac{\sum_{j=1}^{M} \kappa_j p_j (1 + t{\cal G}_j \widetilde{k}_j)}{\sum_{j=1}^M p_j (1 + t{\cal G}_j \widetilde{k}_j)} \leq \kappa_l.\end{array}\right. \end{equation}

\noindent Note that, for fixed $t$,
${\cal G} : \mathbb{S}^M \rightarrow \{0\} \times \mathbb{R}^{M-2}_+ \times \{0\}$, where 
\begin{equation}\label{eqesem} \mathbb{S}^M = \left \{ (p_1, \ldots, p_M) : 0 < p_j < 1, \; j=1, \ldots, M,\; \sum_{j=1}^M p_j = 1 \right \}.
\end{equation}

\paragraph{Notation} Throughout, a discrete target probability $\underline{p} = (p_1, \ldots, p_M)$ will be taken as a {\em row} vector.\vspace{5mm}

\noindent The following lemma shows that such a function is well defined, which is a necessary step in accommodating the killing field. 

\begin{Lmm} \label{lmmgeee} For a given $\underline{p} \in \mathbb{S}^M$,  there exists $({\cal G}_1,\ldots, {\cal G}_{M})$ satisfying Equation~\eqref{eqgeee}.
\end{Lmm}

\paragraph{Proof} 
Consider $(\alpha_1, \ldots, \alpha_M)$ such that $\alpha_j > 0$ for all $j$ and $\sum_{j=1}^M \alpha_j = 1$. Now consider, for some $k \in \{2, \ldots, M-1\}$, $\beta_k \in (0,1)$ and, for $j \neq k$, $\beta_j = \frac{1 - \beta_k}{1 - \alpha_k}\alpha_j$. Then $\sum_{j=1}^M \beta_j =  1$. Let $y = \sum_{j=1}^M \kappa_j \alpha_j$ and let $z = \sum_{j=1}^M \kappa_j \beta_j$, then 

\[ z = y \frac{1 - \beta_k}{1 - \alpha_k} + \frac{(\beta_k - \alpha_k)}{1 - \alpha_k} \kappa_k \Rightarrow y = z + \frac{(\beta_k - \alpha_k)}{(1 - \beta_k)}(z -   \kappa_k).\]

\noindent so that, for $z < \kappa_k$, it follows that $\beta_k < \alpha_k \Rightarrow y < \kappa_k$ and $\beta_k > \alpha_k \Rightarrow y < z < \kappa_k$. It therefore follows that $z < \kappa_k \Rightarrow y < \kappa_k$. 
\vspace{5mm}

\noindent Recall that $\widetilde{k}_1 = \widetilde{k}_M = 0$ and let ${\cal G}_1^+ = {\cal G}_M^+  = {\cal G}_1^- = {\cal G}_M^- = 0$. For $j = 2, \ldots, M-1$ define ${\cal G}_j^+$ and ${\cal G}_j^-$ by: 

\[\left\{ \begin{array}{ll} {\cal G}_j^- =  \frac{1}{tp_j} \frac{(\kappa_{j+1} - \kappa_{j-1})}{(\kappa_{j+1} - \kappa_j)(\kappa_j - \kappa_{j-1})} \sum_{a = 1}^{j-1} (\kappa_j - \kappa_a)p_a (1 + t{\cal G}_a^- \widetilde{k}_a) & 2 \leq j \leq M-1    \\
{\cal G}_j^+ = \frac{1}{tp_j} \frac{(\kappa_{j+1} - \kappa_{j-1})}{(\kappa_{j+1} - \kappa_j)(\kappa_j - \kappa_{j-1})} \sum_{a = j+1}^M (\kappa_a - \kappa_j) p_a (1 + t {\cal G}_a^+ \widetilde{k}_a) & 2 \leq j \leq M-1 \end{array} \right. \] 

\noindent Then these are well defined and positive. Define $x_{0,j}$ by:
\[ \left\{ \begin{array}{ll} x_{0,1} = \frac{\sum_{a=1}^M \kappa_a p_a (1 + t{\cal G}_a^+ \widetilde{k}_a)}{\sum_{a=1}^M p_a (1 + t{\cal G}_a^+ \widetilde{k}_a) } & \\ 
 x_{0,j} = \frac{\sum_{a=1}^{j-1} \kappa_a p_a (1 + t{\cal G}_a^- \widetilde{k}_a) + \sum_{a=j}^M \kappa_a p_a (1 + t{\cal G}_a^+ \widetilde{k}_a)}{\sum_{a=1}^j p_a (1 + t{\cal G}_a^- \widetilde{k}_a) + \sum_{a=j+1}^M p_a(1 + t{\cal G}_a^+\widetilde{k}_a)} & j = 2,\ldots, M \\
x_{0,M+1} = \frac{\sum_{a=1}^M \kappa_a p_a (1 + t{\cal G}_a^- \widetilde{k}_a)}{\sum_{a=1}^M p_a (1 + t{\cal G}_a^- \widetilde{k}_a) }\end{array}\right. \]

\noindent Clearly $\kappa_1 < x_{0,j} < \kappa_M$ for each $j \in \{1, \ldots, M+1\}$. To prove the lemma, it is necessary and sufficient to show that there is a $j$ such that $\kappa_{j-1} < x_{0,j} \leq \kappa_j$.  Note that $x_{0,M} < \kappa_M$. Suppose that $x_{0,j} < \kappa_{j-1}$. Then, it follows from the argument above that  $x_{0,j-1} < \kappa_{j-1}$. If, furthermore, $x_{0,j-1} > \kappa_{j-2}$, then existence has been established; otherwise, proceed inductively. Since $\kappa_1 < x_{0,j} < \kappa_M$ for all $j$, the result follows. 
\qed \vspace{5mm}

\noindent Let 

\begin{equation}\label{eqx02} x_0 = \frac{\sum_{j=1}^{M} \kappa_j p_j (1 + t{\cal G}_j (t,\underline{p}) \widetilde{k}_j)}{\sum_{j=1}^M p_j (1 + t{\cal G}_j(t, \underline{p}) \widetilde{k}_j)}.
\end{equation}

\noindent  This will give the initial condition for the process for geometric / exponential stopping times. $x_0$ may be considered as the average of $\underline{\kappa}$ under the measure $Q(t,\underline{p}, {\cal G}(t,\underline{p}))$ where the quantity $Q$ is defined by~\eqref{eqqespj2} below:

\begin{equation}\label{eqqespj2} Q_j(s,\underline{p},  \underline{\lambda}) = \frac{p_j (1 + s \widetilde{k}_j \lambda_j)}{\sum_{i=1}^M p_i (1 + s \widetilde{k}_i \lambda_i)} \qquad j \in \{1, \ldots, M\}
\end{equation}

\noindent (the killing field is considered fixed; this quantity will be considered as a function of time variable, the target probability and the intensities when it is used later). 

\section{Results} \label{secres} This section states the main results of the article, which are given as  Theorems~\ref{thgeexp2},~\ref{thnbt2},~\ref{thctlim2} and~\ref{thdrk}. These theorems are stated separately, because each of them is of use in its own right. Firstly, Theorem~\ref{thgeexp2} concerns Exponential and Geometric times. In this setting, an explicit solution can be obtained; $\alpha$ and $\underline{\lambda}$ are determined uniquely and there are equations to produce the explicit values. Theorem~\ref{thnbt2} considers Negative Binomial  and Gamma times. Uniqueness is not shown, but the result comes in terms of the solution to an explicit fixed point problem. Theorem~\ref{thctlim2} takes an appropriate limit to obtain the result for deterministic times. While the result of Theorem~\ref{thctlim2} is the objective, the result of Theorem~\ref{thnbt2} which is a step along the way has an interesting interpretation in terms of the Local Variance Gamma Model of Carr~\cite{Carr} and is therefore stated as a theorem in its own right.  

Theorems~\ref{thgeexp2},~\ref{thnbt2} and~\ref{thctlim2} consider discrete state spaces, while Theorem~\ref{thdrk} considers arbitrary probability measures where the measure and drift satisfy Hypothesis~\ref{hybmu}.  \vspace{5mm}

\noindent For Theorems~\ref{thgeexp2},~\ref{thnbt2} and~\ref{thctlim2},  let $\underline{p} = (p_1, \ldots, p_M) \in \mathbb{S}^M$, defined by~\eqref{eqesem}. Let ${\cal S} = \{i_1, \ldots, i_M\}$, $i_1 < \ldots < i_M$ be a finite state space with $M$ elements and let $\underline{b} = (b_2, \ldots, b_{M-1})$ satisfy Hypothesis~\ref{hyb}. Let $\underline{k} = (k_2,\ldots, k_{M-1}) \in  \mathbb{R}_+^{M-2}$ denote the killing field and let $\widetilde{\underline{k}} = (\widetilde{k}_1, \ldots, \widetilde{k}_{M})$ be defined by~\eqref{eqtildek}.

\begin{Th} \label{thgeexp2}    There exists a unique $\underline{\lambda} = (\lambda_2, \ldots, \lambda_{M-1}) \in \mathbb{R}^{M-2}_+$,  $\alpha \in (0,1)$, $l \in \{2, \ldots, M\}$ and  $\beta \in (0,1]$ such that for all

\[ h  \in \left ( 0, \min_{j \in \{2, \ldots, M-1\}} \left(\frac{1}{\lambda_j (1 + \widetilde{k}_j) }\right) \right ] \]

\noindent $\widetilde{P}^{(h)}$ given by Definition~\ref{deftrmat} is the one step transition matrix   for a Markov chain $X^{(h)}$ with state space ${\cal S} = \{i_1, \ldots, i_M, D\}$ and time step length $h$ which  satisfies

\[ \alpha p_j = (1 - \beta) \mathbb{P}\left (X_{h\tau}^{(h)} = i_j | X_0^{(h)} = i_{l-1} \right ) + \beta \mathbb{P} \left (X_{h\tau}^{(h)}= i_j | X_0^{(h)} = i_l \right ) \qquad j = 1,\ldots, M, \]

\noindent where $\tau$ is independent of $X^{(h)}$ and satisfies $\tau \sim Ge(a)$ with $a = \frac{t}{t + h}$, so that $\mathbb{E}[h \tau] = \frac{a h}{1-a} = t$. The constant $\alpha$ satisfies:

\begin{equation}\label{eqalph2} \alpha = \frac{1}{1 + t\sum_{j=1}^M p_j \lambda_j  \widetilde{k}_j}\end{equation}

\noindent where $\lambda_1 = \lambda_M = 0$, while $\beta$ satisfies:

\begin{equation}\label{eqbet2} \beta = \frac{x_0 - \kappa_{l-1}}{\kappa_l - \kappa_{l-1}} 
\end{equation}

\noindent where $x_0$ is defined by~\eqref{eqx02}. The intensity vector $\underline{\lambda}$ satisfies:

\begin{equation}\label{eqlamexpsol2} \lambda_j = {\cal G}_j \qquad j = 2, \ldots, M-1 \end{equation}

\noindent where ${\cal G}$ is defined by~\eqref{eqgeee}, existence of such a function given by Lemma~\ref{lmmgeee}. Taking $h \rightarrow 0$, there exists a continuous time Markov chain $Y$ with state space ${\cal S}$, where for each $j = 2, \ldots, M - 1$, site $i_j$ has holding intensity $\lambda_j$ given by the same formula, and $Y$ satisfies

\[ \alpha p_j = (1-\beta) \mathbb{P}(Y_\tau = i_j | Y_0 = i_{l-1}) + \beta \mathbb{P}(Y_\tau = i_j | Y_0 = i_l), \]

\noindent $\tau \sim Exp(1/t)$ (that is, exponential, with expected value $\mathbb{E}[\tau] = t$), $\alpha$ satisfies~\eqref{eqalph2} and $\beta$ satisfies~\eqref{eqbet2}. Let $\underline{a}$ satisfy~\eqref{eqaj} with $\underline{\lambda}$ given by~\eqref{eqlamexpsol2}, then the infinitesimal generator of the process $Y$ is given by: 

\begin{equation}\label{eqinfgen2} a \left(\frac{1}{2} \Delta + b \nabla  - k\right).\end{equation}
\end{Th}

\noindent The quantity $\beta$ is interpreted in the following way: the process $X^{(h)}$ has initial condition $X^{(h)}_0 = x_0 \in (i_{l-1}, i_l]$ such that

\[ \mathbb{P}(X_{0+}^{(h)} = i_l | X_0^{(h)} = x_0) = \beta \qquad \mathbb{P}(X_{0+}^{(h)} = i_{l-1}|X_0^{(h)} = x_0) = (1 - \beta).\]

\noindent Now consider negative binomial times.
 
 \begin{Th}\label{thnbt2} For any $r \geq 1$, there exists an $l \in \{2, \ldots, M\}$,  a vector $\underline{\lambda} = (\lambda_2, \ldots, \lambda_{M-1}) \in \mathbb{R}_+^{M-2}$,  a $\beta \in (0,1]$,  and an $h_0 \in (0,1)$, such that for all $h \in (0,h_0)$, there is an $\alpha$ satisfying 
 \begin{equation}\label{eqalph} \alpha \in \left [ \left( 1 + \frac{(k\lambda)^* t}{r}\right)^{-r} ,1 \right ]\end{equation}
 
 \noindent where  $(k\lambda)^* = \max_{j \in \{2, \ldots, M-1\}} \widetilde{k}_j \lambda_j$ and    $\widetilde{P}^{(h)}$ from Definition~\ref{deftrmat} is the one step transition matrix for  a time homogeneous discrete time Markov chain $X^{(h)}$, time step length $h$ such that

\[ (1 - \beta) \mathbb{P}(X_{h\tau}^{(h)} = i_j | X_0^{(h)} = i_{l-1}) + \beta \mathbb{P}(X_{h\tau}^{(h)} = i_j | X_0^{(h)} = i_l) = \alpha p_j \qquad j = 1,\ldots, M\] where $\tau \sim NB(r, \frac{t}{t + hr})$, so that $\mathbb{E} \left[ h \tau \right ] = t$.

By taking the limit $h \rightarrow 0$, there is a   continuous time, time homogeneous Markov chain $Y$ with transition intensity matrix $\Theta$ given by~\eqref{eqthetmat}, Definition~\ref{defintmat}, such that

\[ (1 - \beta) \mathbb{P}(Y_T = i_j | Y = i_{l-1}) +  \beta \mathbb{P}(Y_T = i_j |Y_0 = i_l) = \alpha p_j \qquad j = 1, \ldots, M\] 

\noindent where $\alpha$ satisfies~\eqref{eqalph} and $T \sim \mbox{Gamma}(r, \frac{r}{t})$; that is, $T$ is a Gamma time, with density function

\begin{equation}\label{eqgamden} f_T(x) = \frac{r^r}{t^r}\frac{1}{\Gamma(r)}x^{r-1}e^{-xr/t} \qquad x \geq 0
\end{equation}

\noindent and expected value  $\mathbb{E}[T] = t$.
 \end{Th}
 
 \noindent This is extended to deterministic time:
 
 \begin{Th}\label{thctlim2}   For a given $t > 0$, there exists a vector $\underline{\lambda} = (\lambda_2, \ldots, \lambda_{M-1}) \in \mathbb{R}_+^{M-2}$, an

\begin{equation}\label{eqalthcl2}\alpha \in \left [  \exp \left\{-t(k \lambda)^*\right\}, 1 \right ],\end{equation}

\noindent where  $(k\lambda)^* = \max_{j \in \{2, \ldots, M-1\}} \lambda_j \widetilde{k}_j$, an $l \in \{2, \ldots, M\}$, a  $\beta \in (0,1]$,   such that  $\Theta$ (Equation~\eqref{eqthetmat} Definition~\ref{defintmat})  is the intensity  matrix for  a time homogeneous continuous time Markov chain $X$ such that

\[ \beta \mathbb{P}(X_t = i_j | X_0 = i_l) + (1 - \beta) \mathbb{P}(X_t = i_j | X_0 = i_{l-1}) = \alpha p_j \qquad j = 1, \ldots, M.\]

\noindent Again, if $\underline{a}$ satisfies~\eqref{eqaj}, then the infintesimal generator of $X$ is defined by~\eqref{eqinfgen2}.
\end{Th}
 
 \noindent Finally,  the continuous limit in the space variable can be taken.
 
\begin{Th}\label{thdrk} Let $\mu$ be a probability measure on $\mathbb{R}$, $b : \mathbb{R} \rightarrow \mathbb{R}$ a drift field and $k$ a killing field such that $(\mu, b, k)$ satisfy  Hypothesis~\ref{hybmu}, then there exists a string measure $m$, a point $x_0 \in \mathbb{R}$, an $\alpha \in (0,1]$ and a function $K$ satisfying Equation~\eqref{eqtildk}  such that $\frac{1}{2}\frac{\partial^2}{\partial m \partial x} + \frac{\partial B}{\partial m}\nabla_m - \frac{\partial K}{\partial m}$ is the infinitesimal generator of a process $X$ which satisfies
\[ \mathbb{P}(X_t \leq x | X_0 = x_0) =  \alpha \mu((-\infty, x]).\]
\end{Th}

\noindent Here $B$ is defined by~\eqref{eqtildb}. The initial condition $X_0 = x_0$ is interpreted as follows: let ${\cal S}$ denote the support of $\mu$. Let $z_- = \sup\{y < z|y \in {\cal S}\}$ and $z_+ = \inf\{y > z | y \in {\cal S}\}$. Then there is a $\beta \in (0,1]$ such that 
\[ \beta \mathbb{P}(X_t \leq x|X_0 = x_{0+}) + (1 - \beta) \mathbb{P}(X_t \leq x | X_0 = x_{0-}) = \alpha \mu((-\infty, x]).\]

\noindent That is, if $x_0 \not \in {\cal S}$, then the process immediately jumps into ${\cal S}$, taking values $x_{0+}$ or $x_{0-}$ with probabilities $\beta$ and $1 - \beta$ respectively:

\[ \mathbb{P}(X_{0+} = x_{0+}|X_0 = x_0) = \beta \qquad \mathbb{P}(X_{0+} = x_{0-}|X_0 = x_0) = 1 - \beta.\]

\noindent Note: if $m$ has a density $m^\prime$, then the infinitesimal generator may be written as

\[ \frac{1}{m^\prime}\left(\frac{1}{2}\frac{\partial^2}{\partial x^2} + b\frac{\partial}{\partial x} - k\right).\]

\section{Proofs of the results in the absence of a killing field}\label{seccoord}

For Theorems~\ref{thgeexp2},~\ref{thnbt2} and~\ref{thctlim2} which consider a finite state space ${\cal S} = \{i_1, \ldots, i_M\}$, let $\underline{\kappa}$ be defined by~\eqref{eqkappa}. For $\underline{\kappa}$ so defined, the quantities $q_{j,j+1}$ and $q_{j,j-1}$ from~\eqref{eqdefq} satisfy~\eqref{eqkap1}.  With $k \equiv 0$, the problem is therefore that of finding a {\em martingale} generalised diffusion when viewed in the changed co-ordinates described above and is therefore solved in the article~\cite{N1}.

For Theorem~\ref{thdrk}, the proof also follows similarly to that of~\cite{N1}, with the following alterations. As in~\cite{N1}, at stage $N$, the points $i_{N,1} < \ldots < i_{N,M_N}$ are chosen in the following way: let

\begin{equation}\label{eqindef}\left\{ \begin{array}{l} i_{N,1} = \inf \left \{x | \mu((-\infty,x]) \geq \frac{1}{2^N} \right \} \\  
i_{N,j} = \inf \left \{x > i_{N,j-1} | \mu((i_{N,j-1},i_{N,j}]) \geq \frac{1}{2^N} \right \} \\  M_N = \inf\left \{j : 1 - \mu((-\infty, i_{N,j}]) \leq \frac{1}{2^N} \right\}. \end{array}\right.
\end{equation}  

\noindent As in~\cite{N1}, let $\underline{p}^{(N)}$ be defined as

\begin{equation}\label{eqpenj} p_j^{(N)} = \left\{ \begin{array}{ll} \mu((-\infty, i_{N,1}]) & j = 1 \\ \mu((i_{N,j-1}, i_{N,j}])  & j = 2, \ldots, M_N - 1 \\
1 - \mu((-\infty, i_{N,M_N}]) & j = M_N.
\end{array} \right.
\end{equation}

\noindent Let $\widetilde{b}$ be defined by~\eqref{eqtildb} and set: 

\begin{equation}\label{eqdd}   b^{(N)}_j = \frac{1}{i_{N,j+1} - i_{N,j-1}}\int_{i_{N,j-1}+}^{i_{N,j+1}-} \widetilde{b}(x) dx \qquad j = 2, \ldots, M_N - 1  \end{equation}

\noindent where $\int_{a+}^{b-}$ means $\int_{(a,b)}$, the integral over the open interval. Note that~\eqref{eqbhyp} of Hypothesis~\ref{hyb} is satisfied if:
\[ \left\{ \begin{array}{l} -1 < \min_{j \in \{2, \ldots, M_{N}-1\}} \frac{i_{N,j+1} - i_{N,j}}{i_{N,j+1} - i_{N,j-1}}\int_{i_{N,{j-1}+}}^{i_{N,j+1}-} \widetilde{b}(x) dx,\\ \max_{j \in \{2, \ldots, M_{N}-1\}} \frac{i_{N,j} - i_{N,j-1}}{i_{N,j+1} - i_{N,j-1}}\int_{i_{N,{j-1}+}}^{i_{N,j+1}-} \widetilde{b} (x) dx < 1. \end{array}\right. \]

\noindent  

\noindent   From~\eqref{eqhypb} and~\eqref{eqhypb2} of Hypothesis~\ref{hybmu}, it follows that there is an $N_0$ such that for all $N > N_0$

\begin{equation}    -1 + \gamma < \min_{j \in \{2, \ldots, M_{N}-1\}}  \int_{i_{N,{j-1}+}}^{i_{N,j+1}-} \widetilde{b}(x) dx \leq \max_{j \in \{2, \ldots, M_{N}-1\}}  \int_{i_{N,{j-1}+}}^{i_{N,j+1}-}  \widetilde{b} (x) dx < 1 - \gamma \label{eqhypb3} \end{equation}

\noindent where $\gamma$ is from~\eqref{eqgamdef}. It follows that ~\eqref{eqbhyp} of Hypothesis~\ref{hyb} is  satisfied for $N > N_0$. For the remainder of the argument, only $N > N_0$ is considered. Using $\underline{b}^{(N)}$ defined by~\eqref{eqdd}, let $\underline{\lambda}_N = (\lambda_{N,2}, \ldots, \lambda_{N,M_N-1})$ ($\lambda_{N,j}$ the holding intensity for site $i_{N,j}$, $j = 2, \ldots, M_N-1$) denote the intensity vector that provides a solution to the marginal distribution problem.  Let

\begin{equation}\label{eqaaa} 
a^{(N)}_j = \left\{ \begin{array}{ll} \lambda_j^{(N)} (i_{N,j+1} - i_{N,j})(i_{N,j} - i_{N,j-1}) & j = 2, \ldots, M_N - 1 \\ 0 & j = 1, M_N \end{array}\right.
\end{equation}

\noindent and

\begin{equation}\label{eqdiscmeas} m^{(N)}(\{i_{N,j}\}) = \left\{ \begin{array}{ll} \frac{(i_{N,j+1} - i_{N,j-1})}{a_j^{(N)}} & j \in \{2, \ldots, M_N - 1\} \\ +\infty &  j = 1, M_N \end{array}\right. 
\end{equation}

\noindent The measure $m^{(N)}$ has support $\{i_{N,1}, \ldots, i_{N,M_N}\}$. Let $z_- = \inf\{x : x \in \mbox{suppt}(\mu)$ and $z_+ = \sup\{x : x \in \mbox{suppt}(\mu)$. Then, from the arguments of~\cite{N1}, there is a limiting measure $m$ such that for any $z_- < x < y < z_+$ there exists a subsequence $(N_j)_{j \geq 1}$ satisfying  

\begin{equation}\label{eqlimmeas} \lim_{j \rightarrow +\infty} \sup_{ x < a < b < y} \left | m^{(N_j)}([a,b]) - m([a,b]) \right | = 0.\end{equation}

\noindent Using the notation of Section~\ref{subcoc}, let $\delta_{N,j} = i_{N,j} - i_{N,j-1}$ for $j = 2, \ldots, M_N$,  $\epsilon_{N,j} = \kappa_{N,j} - \kappa_{N,j-1}$ and  

\begin{equation}\label{eqemminemplusN}
e_{N-} = \inf \left \{j : \sum_{i=1}^j p_i^{(N)} \geq \alpha \right \}, \qquad e_{N+} = \sup \left \{j : \sum_{i= j}^{M_N} p_i^{(N)} \geq \alpha \right \}.
\end{equation} 

\noindent where $0 < \alpha < 0.5$ is a number chosen such that there exists an $a$ and a $b$ such that $a < b$ and $\alpha \geq  \mu((-\infty,a])$ and $\alpha \geq \mu([b,+\infty))$ and an $N_0$ such that 
\begin{equation}\label{eqennought} \inf_{N > N_0}|i_{N,e_{N+}} - i_{N,e_{N-}}| > 0
\end{equation}

\noindent (strict inequality). Only $N > N_0$ where this condition and~\eqref{eqbhyp} are satisfied will be considered. Let 

\begin{equation}\label{eqken} K_N = i_{N,e_{N+}} - i_{N,e_{N-}}
\end{equation}

\noindent and 

\begin{equation}\label{eqeps2N} \epsilon_{N,j+1} = K_N 
\frac{\delta_{N,j+1} \prod_{k=1}^j \left(\frac{1 - \delta_{N,k} b_k^{(N)} }{1 + \delta_{N,k+1} b_k^{(N)}}\right)}
{\sum_{a=e_{N-}}^{e_{N+}-1} \delta_{N,a+1} \prod_{k=1}^a \left(\frac{1 - \delta_{N,k} b_k^{(N)}}{1 + \delta_{N,k+1} b_k^{(N)}}\right ) }.
\end{equation}

\noindent Let $e_N$ satisfy: $\sum_{i=1}^{e_N} p_j^{(N)} \geq \frac{1}{2}$ and $\sum_{i= e_N}^{M_N} p_j^{(N)} \geq \frac{1}{2}$. Set  

\begin{equation}\label{eqkappaN} 
\left\{ \begin{array}{ll} \kappa_{N,e_N} = 0 & \\  \kappa_{N,e_N +j + 1} = \kappa_{N,e_N + j} + \epsilon_{N, e_N + j+1}  & j = 0, \ldots, M_N -e_N -1 \\ \kappa_{N,e_N - j -1} = \kappa_{N,e_N - j} - \epsilon_{N,e_N - j} & j = 0, \ldots, e_N - 2. 
\end{array}\right. 
\end{equation}

\noindent Now note that 
\[ \left\{ \begin{array}{l} i_{N,e_{N+}} \rightarrow c_+ := \sup\left \{x \in \overline{\mbox{suppt}(\mu)} : \mu((-\infty , x)) < 1- \alpha \right \} \\  i_{N,e_{N-}} \rightarrow c_- := \inf \left \{x \in \overline{\mbox{suppt}(\mu)} : \mu((-\infty,x]) \geq \alpha \right \} \end{array}\right. \] 

\noindent so that $K_N \rightarrow c_+ - c_-$, a well defined positive limit and that, by construction, $\kappa_{N,e_{N+}} - \kappa_{N,e_{N-}} = K_N$ for each $N$.\vspace{5mm}

\noindent The function of~\eqref{eqhypb} of Hypothesis~\ref{hybmu}  is to ensure that in the new co-ordinates, the process has a well defined expected value. The following lemma demonstrates that the hypothesis is sufficient for this purpose.  

\begin{Lmm}\label{lmmkpbd} With $b^{(N)}_j$ defined by~\eqref{eqdd} and $(\mu, b)$ satisfying~\eqref{eqhypb}, there exists an $N_0 \in \mathbb{Z}_+$ such that 

\begin{equation}\label{eqkpbd} \sup_{N \geq N_0} \sum_{j=1}^{M_N} |\kappa_{N,j}|p^{(N)}_j < +\infty.
\end{equation}
\end{Lmm}

\paragraph{Proof of Lemma~\ref{lmmkpbd}} Let $\gamma$ be defined by~\eqref{eqgamdef}. Let $N_0$ satisfy: for all $N > N_0$, both~\eqref{eqhypb3} and~\eqref{eqennought} hold where $e_{N+}$ and $e_{N-}$ denote the indices defined in (\ref{eqemminemplusN}). Let 

\[ C_N = K_N \frac{\prod_{k=1}^{e_N  - 1} \left(\frac{1 - \delta_{N,k} b_k^{(N)}}{1 + \delta_{N,k+1} b_k^{(N)}}\right)}{\sum_{a=e_{N-}}^{e_{N+} - 1} \delta_{N,a+1} \prod_{k=1}^{a} \left(\frac{1 - \delta_{N,k} b_k^{(N)}}{1 + \delta_{N,k+1} b_k^{(N)}}\right)}. \]

\noindent where, as above, $K_N = i_{N,e_{N+}} - i_{N,e_{N-}} = \kappa_{N,e_{N+}} - \kappa_{N,e_{N-}}$. Recall the definition of $\kappa_{N,.}$ given by~\eqref{eqkappaN}, that $e_N$ is the index such that $\kappa_{N,e_N} = 0$. Also, $\delta_{N,j} = (i_{N,j} - i_{N,j-1})$. Recall the definition of $b^{(N)}$ from~\eqref{eqdd}. Then for $j > e_N$,

\begin{eqnarray} \kappa_{N,j} &=& \nonumber C_N \sum_{k= e_N}^{j-1} (i_{N,k+1} - i_{N,k})\prod_{l = e_N}^k \left(\frac{1 - \delta_{N,l} b_l^{(N)}}{1 + \delta_{N,l+1} b_l^{(N)}}\right)\\
 & = &  C_N \sum_{k= e_N}^{j-1} (i_{N,k+1} - i_{N,k})\exp \left\{ \sum_{l = e_N}^k \ln \left (1 - \delta_{N,l} b_l^{(N)} \right ) - \ln \left (1 + \delta_{N,l+1} b_l^{(N)} \right )\right\}. \label{eqkapenbd}
\end{eqnarray}

\noindent Similarly, for $j < e_N$, so that $\kappa_{N,j} < 0$,

\begin{equation} - \kappa_{N,j} = C_N\sum_{i = j}^{e_N - 1} (i_{N,k+1} - i_{N,k})\exp \left\{ \sum_{l = i}^{e_N-1} \ln \left (1 - \delta_{N,l} b_l^{(N)} \right ) - \ln \left (1 + \delta_{N,l+1} b_l^{(N)} \right )\right\}. \label{eqkapenbd2}\end{equation}

\noindent Note that 

\begin{eqnarray*} C_N &=& K_N \left(\sum_{a = e_{N-}}^{e_N - 1} \delta_{N,a+1} \prod_{k= a+1}^{e_N-1} \left(\frac{1 + \delta_{N,k+1} b_k^{(N)}}{1 - \delta_{N,k} b_k^{(N)}}\right) + \sum_{a= e_N}^{e_{N_+ - 1}} \delta_{N,a+1} \prod_{k=e_N}^{a} \left( \frac{1 - \delta_{N,k} b_k^{(N)}}{1 + \delta_{N,k+1} b_k^{(N)}} \right)\right)^{-1}\\
&\leq & K_N \left(\sum_{a = e_{N-}}^{e_N - 1} \delta_{N,a+1} \prod_{k= a+1}^{e_N-1} \left(\frac{1 - \delta_{N,k+1} | b_k^{(N)}|}{1 + \delta_{N,k} |      b_k^{(N)}|}\right) + \sum_{a= e_N}^{e_{N_+ - 1}} \delta_{N,a+1} \prod_{k=e_N}^{a} \left( \frac{1 - \delta_{N,k} |b_k^{(N)}|}{1 + \delta_{N,k+1} |b_k^{(N)}|} \right)\right)^{-1}\\
& \leq & \prod_{k=e_{N-}+1}^{e_N-1} \left(\frac{1 + \delta_{N,k} |b_k^{(N)}|}{1 - \delta_{N,k+1}|b_k^{(N)}|}\right) \prod_{k=e_N}^{e_{N+}-1} \left(\frac{1 + \delta_{N,k+1}|b_k^{(N)}|}{1 - \delta_{N,k} |b_k^{(N)}|}\right)
\end{eqnarray*}

\noindent where the equality

\[ \sum_{a=e_{N-}}^{e_{N+}-1} \delta_{N,a+1} = \sum_{a=e_{N-}}^{e_{N+}-1} (i_{N,a+1} - i_{N,a}) = i_{N,e_{N+}} - i_{N,e_{N-}} = K_N\]

\noindent has been used. It follows from~\eqref{eqhypb3}   together with the definition of $c(\gamma)$ given by~\eqref{eqcgamm}, the definition of $b^{(N)}_k$ given by~\eqref{eqdd} and Remarks~\ref{rq1} and~\ref{rq2} about the function $c$ following Hypothesis~\ref{hybmu} that:

\[ C_N \leq  \exp\left\{ 2 \int_{i_{N,e_{N-}}}^{i_{N,e_{N+}}} |\widetilde{b}(x)|dx + 2 c(\gamma)\sum_{k=e_{N-}}^{e_{N+}} \left(\int_{i_{N,k}}^{i_{N,k+1}}  |\widetilde{b}(x) | dx\right)^2\right\}.\]

\noindent  Using $i_{N,e_{N-}} \downarrow c_-$ and $i_{N,e_{N+}} \uparrow c_+$ together with~\eqref{eqhypb} gives  
that $C_N$ is uniformly bounded by a constant $C < +\infty$. It follows that:  

\begin{eqnarray*} \sum_j |\kappa_{N,j}| p_j^{(N)} &\leq& C \left ( \sum_{j= e_N+1}^{M_N} p_j^{(N)}   \sum_{k=e_N}^{j-1} (i_{N,k+1} - i_{N,k})    \right.  \\&& \times \exp\left\{ 2\int_{e_N}^{i_{N,k+1}} |\widetilde{b}(x)|dx  + 2 c(\gamma) \sum_{a=e_N}^k \left(\int_{i_{N,a}}^{i_{N,a+1}} |\widetilde{b}(x)| dx \right)^2 \right\} \\&&
+  \sum_{j=1}^{e_N-1} p_j^{(N)} \sum_{k=j+1}^{e_N} (i_{N,k} - i_{N,k-1})\\&& \left. \hspace{5mm} \times \exp\left\{ 2\int_{i_{N,k-1}}^{e_N} |\widetilde{b}(x)| dx  + 2 c(\gamma)\sum_{a=k-1}^{e_N} \left(\int_{i_{N,a}}^{i_{N,a+1} } |\widetilde{b}(x) | dx \right)^2 \right\} \right ) \\
&& \leq C \int_{-\infty}^\infty \left ( \int_{0 \wedge x}^{0\vee x} e^{F(b,y)} dy\right) \mu(dx)
\end{eqnarray*}

\noindent where $F$ is defined by~\eqref{eqefffF}. The result now follows directly from~\eqref{eqhypb}. \qed \vspace{5mm}

\begin{Lmm}\label{lmmkappdef}  Let $\kappa^{(N)}$ denote the function 

\begin{equation}\label{eqkapen}
 \kappa^{(N)}(x) = \left\{ \begin{array}{lll} \kappa_{N,1} & x < i_{N,1} & \\ \kappa_{N,j} + (\kappa_{N,j+1} - \kappa_{N,j}) \frac{x - i_{N,j}}{i_{N,j+1} - i_{N,j}} & x \in [i_{N,j}, i_{N,j+1}) & 1 \leq j \leq M_{N} - 1 \\ \kappa_{N,M_N} & x \geq i_{N,M_N} & \end{array}\right.
\end{equation} There is a non-decreasing map $\kappa : \mathbb{R} \rightarrow \mathbb{R}$ such that for any $-\infty < a < b < +\infty$,  

\begin{equation}\label{eqkap2}  \lim_{N \rightarrow +\infty} \sup_{x \in [i_{N,1} \vee a, i_{N,M_N} \wedge b]} | \kappa^{(N)}(x) - \kappa(x)| = 0
\end{equation}
\end{Lmm}

\paragraph{Sketch of Proof} Firstly, note that 
\begin{equation}\label{eqderkap} \frac{d \kappa^{(N)}}{dx} = \frac{\kappa_{N,j+1} - \kappa_{N,j}}{i_{N,j+1} - i_{N,j}} = \frac{\epsilon_{N,j+1}}{\delta_{N,j+1}}\qquad x \in [i_{N,j}, i_{N,j+1}).
\end{equation}

\noindent The following argument shows that $\frac{d \kappa^{(N)}}{dx}$ has a well defined limit. From the definition of $\epsilon_{N,.}$, it follows that:

\[ \left. \begin{array}{l} \frac{d\kappa^{(N)}}{dx}(x) =  K_N  \frac{1}{\left ( \sum_{a=j+1}^{e_{N+}-1} \delta_{N,a+1} \prod_{k=j+1}^a \left(\frac{1 - \delta_{N,k} b_k^{(N)}}{1 + \delta_{N,k+1} b_k^{(N)}}\right) + \sum_{a=e_{N-}}^j \delta_{N,a+1} \prod_{k=a+1}^j \left(\frac{1 + \delta_{N,k+1} b_k^{(N)}}{1 - \delta_{N,k} b_k^{(N)}} \right) \right)} \\ x \in [i_{N,j}, i_{N,j+1}) \end{array} \right. \]

\noindent Firstly, $K_N$ has a well defined limit, which is $c_+ - c_- =: K$. Now let ${\cal D}$ denote the set of atoms of $b$. By Hypothesis~\ref{hybmu}, this is countable. For $z \in {\cal D}$, let $z_+ = \inf\{y > z | y \in {\cal D}\}$ and let $z_- = \sup\{y < z | y \in {\cal D}\}$. Let  

\[ j_N(x) = \{ j : x \in [i_{N,j},i_{N,j +1})\}\] and, for $z \in {\cal D}$, let $\Delta \widetilde{b}(z) = \lim_{h \rightarrow 0}\int_{l_-(z) - h}^{l_+(z) + h} \widetilde{b}(x) dx$, where $l_-(z)$ and $l_+(z)$ are defined in the lines above~\eqref{eqhypb2}. Let $\widetilde{b}^c$ denote the continuous part of $\widetilde{b}$ (the part remaining after removing the atoms).  Then, for fixed $x$,

\begin{eqnarray*} \lim_{N \rightarrow +\infty} \frac{d\kappa^{(N)}}{dx}(x) &=& K \left(\int_{c_-}^{x}e^{\int_y^x \widetilde{b}^{(c)}(z) dz}\prod_{z \in {\cal D}: y \leq z \leq x} \left ( \frac{1 + \left(\frac{z_+ - z}{z_+ - z_-} \right) \Delta \widetilde{b}(z)}{1 -  \left( \frac{z - z_-}{z_+ - z_-} \right) \Delta \widetilde{b}(z)}\right)dy \right. \\&& \left. + \int_{x}^{c_+} e^{-\int_x^y \widetilde{b}^{(c)}(z)dz} \prod_{z \in {\cal D}: x \leq z \leq y} \left(\frac{1 - \left(\frac{z - z_-}{z_+ - z_-}\right)\Delta \widetilde{b}(z)}{1 + \left(\frac{z_+ - z}{z_+ - z_-}\right)\Delta \widetilde{b}(z)}\right) dy \right)^{-1}
\end{eqnarray*}

\noindent Let $z_- = \inf\{x \in \mbox{suppt}(\mu)\}$ and $z_+ = \sup\{x \in \mbox{suppt}(\mu)\}$. Conditions~\ref{hycond2} and~\ref{hycond3} of Hypothesis~\ref{hybmu} ensure that this limit is well defined on $[z_-,z_+]$. Furthermore, $i_{N,e_N}$ has a well defined limit and $\kappa^{(N)}(i_{N,e_N}) = 0$ for each $N$. The result now follows almost directly.   
\qed \vspace{5mm}

\noindent Let $X^{(N)}$ denote the process generated by $a^{(N)}\left(\frac{1}{2} \Delta_{{\cal S}_N} + b^{(N)}\nabla_{{\cal S}_N} \right )$ where $\Delta_{{\cal S}_N}$ and $\nabla_{{\cal S}_N}$ are the Laplacian and gradient operators respectively defined on ${\cal S}_N = \{i_{N,1}, \ldots, i_{N,M_N}\}$ (Definition~\ref{deflapder}) and, with $m$ satisfying~\eqref{eqlimmeas}, let $X$ denote the process generated by $\left(\frac{1}{2}\frac{\partial^2}{\partial m \partial x} + \frac{\partial B}{\partial m} \nabla_m\right)$. Let $Y^{(N)} = \kappa^{(N)}(X^{(N)})$ where $\kappa^{(N)}$ is defined by~\eqref{eqkapen} and the mapping $\kappa$ by~\eqref{eqkap2} and let $Y = \kappa (X)$. 
Then $Y^{(N)}$ is a process with state space ${\cal R}_N = \{\kappa_{N,1}, \ldots, \kappa_{N, M_N}\}$ where site $\kappa_{N,j}$ has holding intensity $\lambda^{(N)}_j$ for $j = 1, \ldots, M_N$ and, when it jumps from $\kappa_{N,j}$ for $j \in \{2, \ldots, M_N-1\}$, it jumps to $\kappa_{N,j+1}$ with probability $\frac{\kappa_{N,j} - \kappa_{N,j-1}}{\kappa_{N,j+1} - \kappa_{N,j-1}}$ and to  $\kappa_{N,j-1}$ with probability $\frac{\kappa_{N,j+1} - \kappa_{N,j}}{\kappa_{N,j+1} - \kappa_{N,j-1}}$.  In short, it is a process with infinitesimal generator $\frac{\widetilde{a}^{(N)}}{2}\Delta_{{\cal R}_N}$, where $\Delta_{{\cal R}_N}$ denotes the Laplace operator defined on functions on ${\cal R}_N$ (Definition~\ref{deflapder}) and (with reduction in the notation which is clear)

\[ \widetilde{a}^{(N)}_j := \widetilde{a}^{(N)}(\kappa_{N,j}) = \left\{ \begin{array}{ll}  \lambda_j^{(N)} (\kappa_{N,j+1} - \kappa_{N,j})(\kappa_{N,j} - \kappa_{N,j-1}) & j = 2, \ldots, M_N - 1 \\ 0 & j = 1, M_N.\end{array}\right. \] 

\noindent Let $\widetilde{m}^{(N)}$ denote the measure supported on $\{\kappa_{N,1}, \ldots, \kappa_{N,M_N}\}$ defined by 

\[\left\{ \begin{array}{l} \widetilde{m}^{(N)}(\{ \kappa_{N,j}\}) = \frac{(\kappa_{N,j+1} - \kappa_{N,j-1})}{\widetilde{a}^{(N)}_j} \\ \widetilde{m}^{(N)}(\{\kappa_{N,1}\}) = \widetilde{m}^{(N)}(\{ \kappa_{N,M_N}\}) = +\infty \end{array}\right. \] 

\noindent It follows from the convergence results of~\eqref{eqlimmeas} and~\eqref{eqkap2}  that there is a limit $\widetilde{m}$ such that for the convergent subsequence of~\eqref{eqlimmeas} 

\begin{equation}\label{eqlimmeas2}
 \lim_{j \rightarrow +\infty} \sup_{ \kappa(x) < a < b < \kappa(y)   } \left | \widetilde{m}^{(N_j)}([a,b]) - \widetilde{m}([a,b]) \right | = 0. 
\end{equation}

\noindent As in~\cite{N1}, convergence of processes is based on the following result, which is stated in Kotani-Watanabe~\cite{KW}: 

\begin{Th}[Characterisation of generalised diffusion]\label{thchgd} Let $W$ denote a standard Wiener process starting from $0$ and let $\phi(s,a)$ denote its local time at site $a \in \mathbb{R}$, at time $s \geq 0$. Let $m$ be a measure on $\mathbb{R}$. Let 

\[ T(z,s) = \int_{\mathbb{R}} \phi (s,a-z) m(da)  \]

\noindent and 
\[T^{ -1}(z,s) = \inf\left\{r | \int_{\mathbb{R}} \phi(r,a-z) m(da) \geq s \right \}.\]

\noindent Then $Y(t,z) = z + W(T^{-1}(z,t))$ is a strong Markov process with infinitesimal generator $\frac{1}{2}\frac{\partial^2}{\partial m \partial x}$.
\end{Th} \qed

\noindent Let $f_0^{(N)} = \sum_{j=1}^{M_N} \kappa_{N,j} p^{(N)}_j$. It follows from Equation (\ref{eqkap2}) and the definition of $\underline{p}^{(N)}$ (Equation (\ref{eqpenj})) that there exists an $f_0$ such that $\lim_{N \rightarrow +\infty} |f_0^{(N)} - f_0| = 0$. It therefore follows from the convergence result \eqref{eqlimmeas2} together with Theorem~\ref{thchgd}, that there is a subsequence such that for all $\epsilon > 0$ 

\[\lim_{j \rightarrow +\infty} \mathbb{P} \left(\sup_{0 \leq s \leq t} \left | Y^{(N_j)}_s(f_0^{(N_j)}) - Y_s(f_0) \right | > \epsilon \right) = 0,\]

\noindent where an initial condition $y$ for $Y_s^{(N)}(y)$ is interpreted as:

\[ \left\{ \begin{array}{l} \mathbb{P}\left (Y_{0+}^{(N)}  = \kappa_{N,l_N(y)} | Y_0^{(N)}  = y \right ) = \frac{y - \kappa_{N,l_N(y)-1}}{\kappa_{N,l_N(y)} - \kappa_{N,l_N(y)-1}} =: \beta_N \\ \mathbb{P} \left (Y_{0+}^{(N)}  = \kappa_{N, l_N(y)-1} | Y_0^{(N)}  = y \right ) = \frac{\kappa_{N, l_N(y)} - y}{\kappa_{N, l_N(y)} - \kappa_{N, l_N(y)-1}} = 1 - \beta_N \end{array}\right. \]
\noindent and $l_N(y)$ is defined as the index such that $\kappa_{N,l_N(y)-1} < y \leq \kappa_{N,l_N(y)}$. Let $y_- = \lim_{N \rightarrow +\infty} \kappa_{N, l_N(f_0^{(N)})-1}$ and let $y_+ = \lim_{N \rightarrow +\infty} \kappa_{N, l_N(f_0^{(N)})}$. Then, when $y_- < y_+$ where the inequality is strict, the initial condition $f_0$ for process $Y(f_0)$ is interpreted as:

\[ \mathbb{P}(Y_{0+} = y_+ | Y_0 = f_0) = \frac{f_0 - y_-}{y_+ - y_-} =: \beta \qquad \mathbb{P}(Y_{0+} = y_- | Y_0 = f_0) = \frac{y_+ - f_0}{y_+ - y_-} = 1 - \beta.\] 

\noindent From this,

\[\lim_{j \rightarrow +\infty} \mathbb{P} \left(\sup_{0 \leq s \leq t} \left | X^{(N_j)}_s - X_s \right | > \epsilon \right) = 0 \]

\noindent where $X^{(N)}$ satisfies:

\[ \mathbb{P} \left (X^{(N)}_{0+} = i_{N,l_N(f_0^{(N)})} \right ) = \beta_N \qquad \mathbb{P} \left (X^{(N)}_{0+} = i_{N,l_N(f_0^{(N)})-1} \right ) = 1 - \beta_N \]

\noindent and $\left (\beta^{(N_j)}, i_{N, l_{N_j} (f_0^{(N_j)}) - 1}, i_{N, l_{N_j} (f_0^{(N_j)}) } \right ) \stackrel{j \rightarrow +\infty}{\longrightarrow} \left (\beta, x_{0-}, x_{0+} \right )$ and $X$ satisfies:

\[ \mathbb{P}(X_{0+} = x_{0+}) = \beta \qquad \mathbb{P}(X_{0+} = x_{0-}) = 1 - \beta. \]

\noindent It follows that 

\[ \lim_{j \rightarrow +\infty} \sup_x \left |\mathbb{P}\left( X_t^{(N_j)} \leq x\right) - \mathbb{P} \left (X_t \leq x \right ) \right | = 0\]

\noindent and hence that 

\[ \mathbb{P} \left (X_t \leq x \right )  = \mu ((-\infty, x]) \]

\noindent where $X$ is a diffusion process with infinitesimal generator $\left(\frac{1}{2}\frac{\partial^2}{\partial m \partial x} + \frac{\partial B}{\partial m} \nabla_m \right)$ as required.  \qed 

\section{Introducing the Killing Field: Preliminary Results}\label{seckf}  Attention is now turned to the problem of introducing a killing field $k$. The following sections prove the theorems of the article stated in Section~\ref{secres}; this section presents preliminary results and notation. 

The transition from finite state space to arbitrary measure on $\mathbb{R}$ follows the same proof as~\cite{N1}, together with the arguments of Section~\ref{seccoord}, with only a few additions. Some discussion is necessary for modifying the proofs of~\cite{N1} so that they can accommodate killing for geometric / exponential times and then to modify the fixed point theorem so that the transition can be made to negative binomial times. Once negative binomial times are accommodated, the limiting arguments to obtain the result for deterministic time are straightforward and the limiting arguments to obtain the result for arbitrary state space follow directly from the analysis of~\cite{N1}. \vspace{5mm}

\noindent Recall the definitions of $\widetilde{P}$ (Equation~\eqref{eqptild} Definition~\ref{deftrmat}). Let 

\begin{equation}\label{eqen} \widetilde{{\cal N}}_t  =   \frac{t+h}{h} \left (I - \frac{t}{t+h} \widetilde{P}^{(h)}  \right )  \end{equation}

\noindent  For the problem without drift or killing, this quantity appeared crucially in establishing the result for geometric times in~\cite{N1}, with $\frac{1}{1-a} = \frac{t+h}{h}$ giving $a = \frac{t}{t+h}$. 

The entries of $\widetilde{{\cal N}}$ may be computed quite easily and are given in~\eqref{eqtilden2} below:

\begin{equation}\label{eqtilden2} \widetilde{{\cal N}}_{t;j,k}(\underline{k}, \underline{\lambda}) = \left\{ \begin{array}{ll} 1 & k = j = M+1 \\ 0 & j = M+1, \quad k \neq M+1 \\ - t \lambda_j \widetilde{k}_j & 2 \leq j \leq M -1, \quad k = M+1 \\
0 & (j,k) = (1,M+1), \qquad (j,k) = (M,M+1) \\ 1 & (j,k) = (1,1), \qquad (j,k) = (M,M) \\  1 + t   \lambda_j ( 1 + \widetilde{k}_j ) & k = j \quad  2 \leq j \leq M-1 \\ -t\lambda_j \left(\frac{\kappa_j - \kappa_{j-1}}{\kappa_{j+1} - \kappa_{j-1}}\right)  & k = j+1, \quad 2 \leq j \leq M-1 \\ -t \lambda_j \left(\frac{\kappa_{j+1} - \kappa_j}{\kappa_{j+1} - \kappa_{j-1}}\right)   & k = j-1, \quad 2 \leq j \leq M-1 \\   0 & j = 1,  k \in \{2, \ldots, M\} \\    0 & j = M, \quad k \in \{1, \ldots, M-1\} \\ 0 & (j,k) \in \{1, \ldots, M\}^2  \quad |j - k| \geq 2  \end{array}\right.
\end{equation}

\noindent where $\widetilde{\underline{k}}$ is defined by~\eqref{eqtildek}. Note that this is independent of $h$. Let ${\cal N}_t$ denote the $M \times M$ matrix such that ${\cal N}_{t;j,k} = \widetilde{{\cal N}}_{t;j,k}$ for $(j,k) \in \{1, \ldots, M\}^2$.  

\paragraph{Note} The notation will be suppressed; $\widetilde{{\cal N}}_t(\underline{k},\underline{\lambda})$ and ${\cal N}_t(\underline{k},\underline{\lambda})$ will be written as $\widetilde{{\cal N}}$ and ${\cal N}$ respectively. Some particular variables ($t$ or $\underline{\lambda}$) may be introduced if they are of particular concern for the point under discussion.\vspace{5mm}

\noindent It is straightforward to compute that $\sum_{k = 1}^{M+1} \widetilde{{\cal N}}_{t;j,k} = 1$, but there does not seem to be a direct method to control the {\em absolute} values of the entries of the matrix. Control is therefore obtained by using the inverse. One result used in the sequel is that for integer $p \geq 1$, all the entries of $\widetilde{{\cal N}}^{-p}$ are non negative, bounded between $0$ and $1$ and that for each $j$, $\sum_{k = 1}^{M+1} (\widetilde{{\cal N}}^{-p}_{t})_{jk} = 1$. This follows from the following representation. 

\begin{Lmm}\label{lmmnprep}  For integer $p \geq 1$, $\widetilde{{\cal N}}^{-p}$   has representation:  

\begin{equation}\label{eqtmtnm1rep} (\widetilde{{\cal N}}^{-p}_t)_{jk} = \mathbb{P}(Z_T = k | Z_0 = j)  \end{equation}

\noindent where $Z$ is a continuous time Markov chain on $\{1, \ldots, M+1\}$, with intensity matrix $\Theta$ (Equation~\eqref{eqthetmat}, Definition~\ref{defintmat}) and $T \sim \mbox{Gamma}(p,\frac{1}{t})$, using the parametrisation of a Gamma distribution from~\eqref{eqgamden} (that is the sum of $p$ independent exponential variables, each with intensity parameter $1/t$).
\end{Lmm}

\paragraph{Proof of lemma~\ref{lmmnprep}} Let $a = \frac{t}{t + h}$. Then, for   $h < \frac{1}{\max_{j \in \{2, \ldots, M-1\}} \lambda_j (1 + \widetilde{k}_j)}$, where $\widetilde{\underline{k}}$ is defined by~\eqref{eqtildek},

\begin{eqnarray*} (\widetilde{{\cal N}}^{-p}_t)_{i,j} &=& \left(\frac{1}{1-a} (I - a\widetilde{P}^{(h)})\right)^{-p}_{i,j} \\
&=& (1-a)^p \sum_{k=0}^\infty \binom{p + k - 1}{k} a^k (\widetilde{P}^{(h)})^k_{i,j} = (1-a)^p \sum_{k=0}^\infty \binom{p + k - 1}{k} a^k \mathbb{P}(Z^{(h)}_{kh} = j | Z^{(h)}_0 = i) \end{eqnarray*}
 
\noindent  where $Z^{(h)}$ is a Markov chain with state space $\{1, \ldots, M+1\}$ and one-step transition matrix $\widetilde{P}^{(h)}$ (where $h$ is the time step length) defined by~\eqref{eqptild} Definition~\ref{deftrmat}. Since

\[ \mathbb{P}(\tau = k) = \binom{p + k - 1}{k} a^k (1-a)^p \qquad k = 0,1,2, \ldots \]

\noindent is the probability mass function of an $NB(p,a)$ random variable, it follows that 
 
\[ (\widetilde{{\cal N}}^{-p})_{i,j} = \mathbb{P}\left(Z_{h\tau}^{(h)} = j|Z_0^{(h)} = i \right ) \qquad \mbox{where} \qquad \tau \sim NB \left (p, \frac{t}{t + h} \right )\]

\noindent so that $\mathbb{E}[h\tau] = \frac{hp(t/t+h)}{(h/t+h)} = pt$. The fact that $h\tau$ converges in distribution to $T \sim \mbox{Gamma} \left (p, \frac{1}{t} \right )$ as $h \rightarrow 0$, $Z^{(h)}$ converges (in the sense of finite dimensional marginals) to the required continuous time Markov chain $Z$ and $Z^{(h)}_{h\tau} \stackrel{h \rightarrow 0}{\longrightarrow_{(d)}} Z_T$ follows the proof found in~\cite{N1}.  \qed\vspace{5mm} 

\noindent The following precautionary lemma is introduced to deal with a problem that does not arise in~\cite{N1}; it is necessary to establish that the Fixed Point Theorem (Theorem~\ref{thfpt}, which is the heart of the proof) does not give a process that is dead with probability $1$ at the terminal time. 

\begin{Lmm}\label{lmmcruclb} 
 For integer $p \geq 1$, there exists a constant $c > 0$, which is independent of $\underline{\lambda}$, such that for all $j \in \{1, \ldots, M\}$, $(\widetilde{\cal N}_t^{-p})_{j,M+1} < 1 - c$. 
\end{Lmm}

\paragraph{Proof} 
Recall the representation of the previous lemma: $(\widetilde{\cal N}_t^{-p})_{j,k} = \mathbb{P}(Z_T = k | Z_0 = j)$. It follows from Equation~\eqref{eqthetmat} that the corresponding embedded discrete time chain has transitions $p_{j,j+1} = \frac{q_{j,j+1}}{1 + \widetilde{k}_j}$, $p_{j,j-1} = \frac{q_{j,j-1}}{1 + \widetilde{k}_j}$, $p_{j,M+1} = \frac{\widetilde{k}_j}{1 + \widetilde{k}_j}$, where $q_{j,j+1}$ and $q_{j,j-1}$ are defined by~\eqref{eqdefq} (expressed as~\eqref{eqkap1} in the drift free coordinates). These transitions do not depend on $\underline{\lambda}$.  If the process reaches site $1$, it remains there; if the process reaches site $M$ it remains there. By considering a lower bound on the probability that the process reaches site $1$, it follows that, for $j \neq M$, 
\[ 1 - \widetilde{{\cal N}}_{j,M+1} \geq \prod_{i=2}^{M-1} \left(\frac{q_{i,i-1}}{1 + \widetilde{k}_i}\right) > 0\]
\noindent as required. This is a lower bound on the probability that the process never reaches the cemetery site $M+1$. \qed \vspace{5mm}

\noindent The following lemma is used in the Fixed Point Theorem, to show that as $\epsilon \rightarrow 0$, the sequence of fixed points for the approximating problems remains bounded.

 \begin{Lmm}\label{lmmpr2} 
 If $\lambda_j \rightarrow +\infty$ then $({\cal N}_t^{-1})_{.j} \rightarrow \underline{0}$ and consequently $({\cal N}_t^{-p})_{.j} \rightarrow \underline{0}$ for any integer $p \geq 1$ where the notation $.j$ denotes the $j$th column of the matrix. 
\end{Lmm}

\paragraph{Proof} Let $\beta_{lj} = ({\cal N}_t^{-1})_{lj}$. Then $0 \leq \beta_{lj} \leq 1$.   $\beta$ satisfies the following system: 

\begin{equation}\label{eqqlamzer2} -q_{l,l-1}\beta_{l-1,j} + \left((1 + k_l) + \frac{1}{t\lambda_l}\right)\beta_{lj} - q_{k,k+1} \beta_{k+1,j} = \left\{ \begin{array}{ll} 0 & k \neq j \\ \frac{1}{t \lambda_j} & l = j \end{array}\right. \end{equation}

\noindent where $q_{l,l-1}, q_{l,l+1}$ for $l = 2, \ldots, M-1$ is defined by~\eqref{eqdefq} and the following definition is used for $l = 1,M$: 

\[\left\{ \begin{array}{ll} q_{l,l-1}\beta_{l-1,j} = 0 & l = 1 \\ q_{l,l+1}\beta_{l+1,j} = 0 & l = M. \end{array}\right. \]  

\noindent From the a priori bounds on $\beta_{lj}$ (namely $0 \leq \beta_{lj} \leq 1$ and $\sum_j \beta_{lj} = 1$ for each $l$) which follow directly from Lemma~\ref{lmmnprep}, it follows from~\eqref{eqqlamzer2} that $\beta_{lj} \stackrel{\lambda_j \rightarrow +\infty}{\longrightarrow} 0$ for all $l = 1, \ldots, M$. 
\qed \vspace{5mm}

\noindent The following lemma is key to proving Theorem~\ref{thgeexp2}, since the proof of Theorem~\ref{thgeexp2} boils down to solving the system of equations defined by~\eqref{eqthg} given below. 

\begin{Lmm}\label{lmmkey2}  
Let $\underline{p}$ be a probability measure over $\{1, \ldots, M\}$. There exists a unique $\alpha \in (0,1]$,  $l \in \{2, \ldots, M\}$, $x_0 \in (\kappa_{l-1}, \kappa_l]$ and $\underline{\lambda} \in \mathbb{R}_+^{M-2}$  satisfying~\eqref{eqthg}:

\begin{equation}\label{eqthg}
 (\widehat{\underline{p}}\widetilde{{\cal N}}_t)_k = \left\{ \begin{array}{ll} \frac{x_0 - \kappa_{l-1}}{\kappa_l - \kappa_{l-1}} & l = k \\ \frac{\kappa_l - x_0}{\kappa_l - \kappa_{l-1}} & k = l-1 \\ 0 &  \mbox{otherwise} \end{array}\right. 
\end{equation}

\noindent where $\widehat{p}_k = \alpha p_k$ for $k = 1,\ldots, M$  and $\widehat{p}_{M+1} = 1 - \alpha$. The solution is the following: $\alpha$ satisfies Equation (\ref{eqalph2}), $x_0$ satisfies Equation (\ref{eqx02}) and $\underline{\lambda} = (\lambda_2, \ldots, \lambda_{M-1})$ satisfies Equation (\ref{eqlamexpsol2}), where ${\cal G}$ is defined by Equation (\ref{eqgeee}) and 
$Q$ is defined by Equation (\ref{eqqespj2}). 
\end{Lmm}

\paragraph{Proof} The equation given by~\eqref{eqthg} for $k = M+1$ is:
\[ (1 - \alpha)  \widetilde{{\cal N}}_{t;M+1,M+1} + \alpha \sum_{j=1}^M p_j \widetilde{{\cal N}}_{t;j,M+1} = 0,\]

\noindent which is: 

\[ 1 - \alpha = \alpha t \sum_{j=1}^M p_j \lambda_j \widetilde{k}_j,\]

\noindent where $(\widetilde{k}_1, \ldots, \widetilde{k}_{M})$ is  defined by~\eqref{eqtildek}.  It follows that $\alpha \in (0,1]$ is required to satisfy:

\[  \alpha = \frac{1}{1 + t\sum_{j=1}^M p_j \lambda_j \widetilde{k}_j} \] 

\noindent so that, if there is a solution, then $\alpha$ is uniquely determined with this value.  Let $x_0 \in \mathbb{R}$, $l \in \{2, \ldots, M\}$ and let $\underline{v}(l, x_0)$ satisfy: 

\[ v_j (x_0,l) = \left\{ \begin{array}{ll} \frac{x_0 - \kappa_{l-1}}{\kappa_l - \kappa_{l-1}} & j = l \\ \frac{\kappa_l - x_0}{\kappa_l - \kappa_{l-1}} & j = l-1 \\ 0 & j \neq l, l-1. \end{array}\right. \] 

\noindent There are $M$ equations involving the $M-2$ unknowns, $\lambda_2, \ldots, \lambda_{M-1}$.  These equations are:

\begin{equation}\label{eqtheq2} \left\{ \begin{array}{ll} \alpha p_1     - t \lambda_2 \alpha p_2 q_{21} = v_1 & \\
  -t \lambda_{j-1} \alpha p_{j-1}q_{j-1,j} + \alpha p_j (1 + t  \lambda_j (1 +  \widetilde{k}_j)) - t \lambda_{j+1} \alpha p_{j+1} q_{j+1,j} = v_j & j = 2, \ldots, M-1 \\
  -t \lambda_{M-1} q_{M-1,M} \alpha p_{M-1} + \alpha p_M  = v_M & 
\end{array} \right. \end{equation}

\noindent Set $\Lambda_j = t\alpha p_j \lambda_j$ and $\widetilde{p}_j = \alpha p_j ( 1 + t  \lambda_j  \widetilde{k}_j)$. Since 

\[ \alpha =  \frac{1}{1 + t \sum_{j=1}^M p_j \lambda_j \widetilde{k}_j} = \frac{1}{\sum_{j=1}^M p_j (1 + t \lambda_j \widetilde{k}_j)},\] 

\noindent it follows from the definition of $Q$ (Equation ~\eqref{eqqespj2}) that $\underline{\widetilde{p}} = Q (t,\underline{p})$ and $\sum_{j=1}^M \widetilde{p}_j = 1$. The system of equations~\eqref{eqtheq2} may be written, with these values, as~\eqref{eqrewrite}:

\begin{equation}\label{eqrewrite}
\left\{\begin{array}{ll} \widetilde{p}_1 - \Lambda_2 q_{21} = v_1 & \\ - \Lambda_{j-1} q_{j-1,j} + (\widetilde{p}_j + \Lambda_j) - \Lambda_{j+1} q_{j+1,j} = v_j & j = 2, \ldots, M-1 \\
       -\Lambda_{M-1}q_{M-1,M} + \widetilde{p}_M = v_M. & 
      \end{array}\right.
\end{equation}

\noindent which is a linear system of $M$ equations with $M-2$ unknowns. To show that it is of rank at most $M-2$:  summing both left hand side and right hand side give $1$ for any choice of $x_0$. \vspace{5mm}

\noindent Also,  

\[ \sum_{j} \kappa_j v_j = \kappa_l \frac{x_0 - \kappa_{l-1}}{\kappa_l - \kappa_{l-1}} + \kappa_{l-1}\frac{\kappa_l - x_0}{\kappa_l - \kappa_{l-1}} = x_0,\]

\noindent  It follows that $x_0$ is required to satisfy 

\[ x_0 = \alpha \sum_{j=1}^M \kappa_j p_j (1 + t \lambda_j \widetilde{k}_j) = \frac{\sum_{j=1}^M \kappa_j p_j (1 + t \lambda_j \widetilde{k}_j)}{\sum_{j=1}^M p_j (1 + t \lambda_j \widetilde{k}_j)}.\]

\noindent It follows that if there is a solution, then $\alpha$, $l$ and $x_0$ are uniquely determined with the values given in the statement of the lemma. \vspace{5mm}

\noindent Since $\sum_{j=1}^M \widetilde{p}_j = 1$, it follows that the system of equations given by~\eqref{eqrewrite} is that studied in~\cite{N1}. From~\cite{N1}, it follows that $\underline{\Lambda}$ satisfies:  

\[ \underline{\Lambda} = {\cal L}(\widetilde{p} ) \]

\noindent where 

\[ {\cal L}_j(\underline{p} ) = \left\{ \begin{array}{ll} \frac{(\kappa_{j+1} - \kappa_{j-1})}{(\kappa_{j+1} - \kappa_j)(\kappa_j - \kappa_{j-1})} \left(\sum_{k=1}^{j-1} (\kappa_j - \kappa_k)p_k \right) & 2 \leq j \leq l-1 \\ \frac{(\kappa_{j+1} - \kappa_{j-1})}{(\kappa_{j+1} - \kappa_j)(\kappa_j - \kappa_{j-1})} \left(\sum_{k=j+1}^M (\kappa_k - \kappa_j) p _k \right) & l \leq j \leq M-1 \\ 0 & j = 1 \quad \mbox{or} \quad M. \end{array}\right. \]

 \noindent Therefore any solution satisfies~\eqref{eqlampr1}:
 
\begin{equation}\label{eqlampr1}
\lambda_j \widetilde{k}_j = \frac{1}{t\alpha p_j}{\cal L}_j(Q (t,\underline{p}) ) = \frac{1 + t \lambda_j \widetilde{k}_j}{t Q (t,\underline{p},j)}{\cal L}_j(Q (t,\underline{p}) ) = \frac{1}{t}(1 + t\lambda_j \widetilde{k}_j){\cal F}_j (Q (t,\underline{p} )) \qquad j = 2, \ldots, M-1
\end{equation}
 
\noindent where ${\cal F}$ is defined (as in~\cite{N1}) by~\eqref{eqeff}:
 
\begin{equation}\label{eqeff}
{\cal F}_j(\underline{p}) = \left\{\begin{array}{ll} \frac{1}{p_j}{\cal L}_j (\underline{p}) & j = 2, \ldots, M-1 \\ 0 & j = 1, M \end{array}\right. 
\end{equation}

 \noindent That is,  $\underline{\lambda}$ is a solution if and only if $\lambda_1 = \lambda_M = 0$ and for $j = 2, \ldots, M-1$,  

\begin{eqnarray} \nonumber 
 \lambda_j &=& \frac{1}{t}(1 + t \lambda_j \widetilde{k}_j){\cal F}_j (Q(t,\underline{p}))\\
  \nonumber &=& \frac{1}{t}(1 + t\lambda_j \widetilde{k}_j) \frac{{\cal L}_j (Q(t,\underline{p}))}{p_j (1 + t\lambda_j \widetilde{k}_j)}\frac{1}{\alpha}\\  \nonumber 
 &=& \frac{1}{t p_j}\frac{(\kappa_{j+1} - \kappa_{j-1})}{(\kappa_{j+1} - \kappa_j)(\kappa_j - \kappa_{j-1})} \times \left\{ \begin{array}{ll} \sum_{i = 1}^{j-1} (\kappa_j - \kappa_i)p_i (1 + t\lambda_i \widetilde{k}_i) & 2 \leq j \leq l-1 \\ \sum_{i = j+1}^M (\kappa_i - \kappa_j) p_i (1 + t \lambda_i \widetilde{k}_i) & l \leq j \leq M-1 \end{array}\right.\\
 &=& {\cal G}_j \label{eqdksol}
\end{eqnarray}

\noindent where ${\cal G}$ is defined by~\eqref{eqgeee}. The function ${\cal G}$ exists by Lemma~\ref{lmmgeee}. The proof of Lemma~\ref{lmmkey2} is complete. \qed

\section{Stopping at Independent Geometric or Exponential Time}

The purpose of this section is to prove Theorem~\ref{thgeexp2}.  

\paragraph{Proof of Theorem~\ref{thgeexp2}} This is equivalent to existence and uniqueness of an $l \in \{2, \ldots, M\}$, $\beta \in (0,1]$, $\alpha \in (0,1]$ and a $\underline{\lambda} \in \mathbb{R}^{M-2}$ such that $\widetilde{P}^{(h)}$ (Definition~\ref{deftrmat}, Equation~\eqref{eqptild}) is the transion matrix for a chain $X^{(h)}$ such that $\widehat{p}$ defined as:

\begin{equation}\label{eqwhp}
 \widehat{p}_j = \left\{ \begin{array}{ll} \alpha p_j & j = 1, \ldots, M \\ 1 - \alpha & j = M+1 \end{array}\right. 
\end{equation}

\noindent satisfies:

\[ \widehat{p}_j    = (1-a)\left((1 - \beta) (I - a\widetilde{P}^{(h)}  )^{-1}_{l,j} +   \beta (I - a\widetilde{P}^{(h)} )^{-1}_{l-1,j}\right) \]

\noindent where $X_0^{(h)} = x_0$ for some $x_0 \in (i_{l-1}, i_l]$ and $\beta \in (0,1]$ is a number such that  
\[ \beta = \mathbb{P}(X_{0+}^{(h)} = i_{l-1} | X_0^{(h)} = x_0) \qquad (1-\beta ) = \mathbb{P}(X_{0+}^{(h)} = i_l | X_0^{(h)} = x_0),\]

\noindent It follows that 

\[ \frac{1}{1-a}\left( \underline{\widehat{p}}(I -a \widetilde{P}^{(h)} \right)_k =   \frac{t+h}{h} \left( \underline{\widehat{p}}(I - \frac{t}{t+h}\widetilde{P}^{(h)} \right)_k = 
\left\{ \begin{array}{ll} 1 - \beta  & l = k \\   \beta & k = l-1 \\ 0 & \mbox{otherwise}, \end{array}\right. \]

\noindent which is equivalent to showing existence of an $\alpha$, $l$,  $\beta$ and $\underline{\lambda}$ such that 

\[ ( \underline{\widehat{p}} {\cal N})_k = \left\{ \begin{array}{ll} (1 - \beta) & l = k \\  \beta & k = l-1 \\ 0 & \mbox{otherwise}. \end{array}\right. \] 

\noindent The result now follows directly from Lemma~\ref{lmmkey2} with $\underline{\kappa} = (\kappa_1, \ldots, \kappa_M)$ the change of coordinates described in Section~\ref{subcoc} and 

\begin{equation}\label{eqlamsol2}
\left\{
\begin{array}{l}  \alpha = \frac{1}{1 + t\sum_{j=1}^M p_j \lambda_j \widetilde{k}_j}, \qquad \beta = \frac{\kappa_l - x_0}{\kappa_l - \kappa_{l-1}} \qquad x_0 = \frac{\sum_{j=1}^M \kappa_j p_j (1 + t \lambda_j \widetilde{k}_j)}{\sum_{j=1}^M p_j (1 + t \lambda_j \widetilde{k}_j)} \\  \lambda_j = {\cal G}_j ( \underline{p}, \underline{k} , \underline{\kappa})\end{array}\right. 
\end{equation}
\noindent where ${\cal G}$ is defined by~\eqref{eqgeee} and $\underline{\widetilde{k}}$ by~\eqref{eqtildek}.\vspace{5mm}
 
\noindent The result now follows for $0 < h < \frac{1}{\max_{j \in \{2, \ldots, M-1\}} \lambda_j(1 + \widetilde{k}_j)}$. The limiting argument to obtain  a continuous time process as  $h \rightarrow 0$, which has the prescribed marginal when stopped at an exponential time is given in~\cite{N1}. \qed \vspace{5mm}

\noindent The case with drift and killing on a finite state space, where the process is stopped at an independent exponential time, has now been solved.

\section{Negative Binomial, Gamma and Deterministic Time} This section is devoted to the proofs of Theorems~\ref{thnbt2}  and~\ref{thctlim2}. They follow the lines of the proofs in~\cite{N1}, with some additional ideas  required to deal with the killing field.

\subsection{Proof of Theorem~\ref{thnbt2}} This follows by appealing to the fixed point theorem, Theorem~\ref{thfpt}. As before, let $\tau \sim NB(r,a)$, with $a = \frac{t}{t+hr}$, so that $\mathbb{E}[\tau ] = \frac{ra}{1 - a} = t$.  Then, with $\widetilde{P}$ defined by Equation~\eqref{eqptild} Definition~\ref{deftrmat},  

\[ \frac{1}{1-a} (I - a\widetilde{P}^{(h)} ) =  \frac{t+hr}{hr} \left( I - \frac{t}{t + hr} \widetilde{P}^{(h)} \right) = \frac{(t/r) + h}{h}\left( I - \frac{(t/r)}{(t/r) + h} \widetilde{P}^{(h)}  \right ) =  \widetilde{{\cal N}}_{t/r} . \]

\noindent If $\tau \sim NB(r,a)$ with $a = \frac{t}{t + hr}$, then $\underline{\lambda}$ provides a solution if and only if there is an $\alpha \in (0,1)$, an $l$ and an $x_0 \in (\kappa_{l-1}, \kappa_l]$ such that 
 
 \[ \alpha p_j =  \frac{\kappa_{l} - x_0}{\kappa_l - \kappa_{l-1}}\mathbb{P}(X_{h\tau}^{(h)} = i_j | X_0^{(h)} = i_{l-1}) + \frac{x_0 - \kappa_{l-1}}{\kappa_l - \kappa_{l-1}} \mathbb{P}(X_{h\tau}^{(h)} = i_j | X_0^{(h)} = i_l) \quad j = 1, \ldots, M\]
 
 \noindent Let $\widehat{p}_j = \alpha p_j$ for $j = 1, \ldots M$ and $\widehat{p}_{M+1} = 1 - \alpha$. Then $(\alpha, x_0, \underline{\lambda})$ provide a solution if and only if 
 
\begin{eqnarray*} \widehat{p}_j &=& \frac{\kappa_l - x_0}{\kappa_l - \kappa_{l-1}}\sum_{k=0}^\infty  \mathbb{P}(X_{hk}^{(h)} = i_j | X_0^{(h)} = i_{l-1})\mathbb{P}(\tau = k) \\ && \hspace{5mm} + \frac{x_0 - \kappa_{l-1}}{\kappa_l - \kappa_{l-1}}\sum_{k=0}^\infty \mathbb{P}(X_{hk}^{(h)} = i_j | X_0^{(h)} = i_{l}) \mathbb{P}(\tau = k) \\
 &=& (1-a)^r \left( \frac{\kappa_l - x_0}{\kappa_l - \kappa_{l-1}} \sum_{k=0}^\infty \binom{k+r - 1}{k} a^k (\widetilde{P}^{(h)k})_{l-1,j} + \frac{x_0  - \kappa_{l-1}}{\kappa_l - \kappa_{l-1}} \sum_{k=0}^\infty \binom{k + r - 1}{k} a^k (\widetilde{P}^{(h)k})_{l,j} \right )\\
 &=& \frac{\kappa_l - x_0}{\kappa_l - \kappa_{l-1}}((I - a\widetilde{P}^{(h)})^{-r})_{l-1,j} + \frac{x_0 - \kappa_{l-l}}{\kappa_l - \kappa_{l-1}} ((I - a\widetilde{P}^{(h)})^{-r})_{l,j}
\end{eqnarray*}

\noindent Let $\underline{\widehat{v}}$ be the $M+1$ vector and $\underline{v}$ the $M$ vector defined by $v_l = \widehat{v}_l = \frac{x_0 - \kappa_{l-1}}{\kappa_l - \kappa_{l-1}}$, $v_{l-1} = \widehat{v}_{l-1} = \frac{\kappa_l - x_0}{\kappa_{l} - \kappa_{l-1}}$ and $v_j = 0$ for $j \neq l-1, l$.  It follows that a solution is provided by any $\alpha \in (0,1)$, $\underline{\lambda}$ and $x_0$ such that 

\[  \underline{\widehat{p}} \widetilde{{\cal N}}_{t/r}^r  = \underline{\widehat{v}} \]

\noindent holds. Let 

\begin{equation}\label{eqqdef2} \underline{q} =   \frac{1}{\sum_{j,k} p_j {\cal N}_{t/r;j,k}^{(r-1)} }\underline{p} {\cal N}_{t/r}^{(r-1)},
\end{equation}

\noindent then $\underline{\lambda}$ provides a solution for all $h \in \left (0, \frac{1}{\max_{j \in \{2, \ldots, M-1\}}  \lambda_j (1 + \widetilde{k}_j)}\right)$ where $\widetilde{\underline{k}}$ is defined by~\eqref{eqtildek}, if and only if there is an $\alpha \in (0,1]$ such that $ \alpha \underline{q} {\cal N}_{t/r}(\underline{\lambda}) = \underline{v}$. It follows from Lemma~\ref{lmmkey2} that $\underline{\lambda}$ is a solution if and only if

\begin{equation}\label{eqsollmm2} \lambda_j =  {\cal G}_j \left (Q\left(\frac{t}{r},\underline{q}\right ) \right ) \qquad j = 1, \ldots, M 
 \end{equation}

\noindent where ${\cal G}$ is defined by~\eqref{eqgeee} and  $Q$ by~\eqref{eqqespj2}. Here 

\[ \alpha = \frac{1}{1 + \frac{t}{r}\sum_{j=1}^M q_j \lambda_j \widetilde{k}_j} \qquad x_0 = \frac{\sum_{j=1}^M \kappa_j q_j (1 + \frac{t}{r} \lambda_j \widetilde{k}_j)}{\sum_{j=1}^M q_j (1 + \frac{t}{r} \lambda_j \widetilde{k}_j)}.\]

\noindent The existence of a $\underline{\lambda}$ satisfying~\eqref{eqsollmm2} follows from Theorem~\ref{thfpt}, which gives existence of a fixed point. For the bounds on $\alpha$, let $\sigma = \inf\{t : X_t \in \{ D\} \}$, then 
\[ \mathbb{P}(\sigma \geq (n+1)h) \geq (1 - h (\lambda k)^*)^n\]
\noindent so that, for $\tau \sim NB\left(r,\frac{t}{t+hr} \right )$ independent of $\sigma$ and using $a = \frac{t}{t+hr}$, for $h < \frac{1}{(\lambda k)^*}$, 

\begin{eqnarray*} \alpha &=& \mathbb{P}(h \tau < \sigma) = \sum_{n=0}^\infty \mathbb{P}(\sigma > nh|\tau = n)\mathbb{P}(\tau = n)\\
 &=& \sum_{n=0}^\infty \mathbb{P}(\sigma \geq (n+1)h)\mathbb{P}(\tau = n) \geq \sum_{n=0}^\infty (1 - h (\lambda k)^*)^{n} \binom{n+r-1}{n}a^n(1-a)^{r}\\
 &=& \frac{ (1 - a)^{r}}{(1 - a(1- h (\lambda k)^*))^{r}}=  \left(1 + \frac{t(\lambda k)^*}{r}\right)^{-r},
\end{eqnarray*}

\noindent  as required. These results hold for all $h \in \left ( 0, \frac{1}{(\lambda k)^*}\right )$ and hence in the continuous time limit as $h \rightarrow 0$. Details of the convergence of finite dimensional marginals are given in~\cite{N1}. 
 \qed 

\subsection{Proof of Theorem~\ref{thctlim2}}\label{mod} This follows almost directly from the proof of Theorem~\ref{thnbt2}; the problem is to show that when the limit is taken, the result is non-trivial. Let $\underline{\lambda}^{(r)}$ denote a solution for the process stopped at an independent $\mbox{Gamma}(r, \frac{r}{t})$ time (parametrisation: the second parameter is an intensity parameter, as with~\eqref{eqgamden}). Let $\mathbb{P}(X_{T_{r}}^{(r)} \in \{D\}) = 1 - \alpha_{r}$ (where $\{D\}$ denotes the `cemetery'; $1 - \alpha_r$ is the probability that the process has been killed by time $T_r$).   Note that  $\lambda_1 = \lambda_M = 0$ (and hence there is no killing of the process once it has reached sites $i_1$ or $i_M$). By the proof of Lemma~\ref{lmmcruclb}, this implies that $\mathbb{P}(X_{T_{r}}^{(r)} \in \{D\}) > \prod_{j=1}^{M-1} \left(\frac{q_{M-j,M-j-1}}{1 + \widetilde{k}_{M-j}}\right)$, where $q$ is defined by~\eqref{eqdefq}. (This lower bound comes from considering the embedded discrete time process; when it jumps from site $j$, it jumps to $j+1$ or $j-1$ or $D$ with probabilities $\frac{q_{j,j+1}}{1 + \widetilde{k}_j}$, $\frac{q_{j,j-1}}{1 + \widetilde{k}_j}$ and $\frac{\widetilde{k}_j}{1 + \widetilde{k}_j}$ respectively. Once it reaches site $1$, it remains there for all time). This lower bound does not depend on $r$. It follows that $\inf_r \alpha_r   > 0$. 

Now suppose that there is a subsequence $r_k$ and a $j \in \{2, \ldots, M-1\}$ such that $\lambda_j^{(r_k)} \rightarrow 0$ then, in the limit, if the process reaches site $i_j$, it remains there for all time, so that either $p_{j+1} = \ldots = p_M = 0$ (if $x_0 \leq i_j$) which is a contradiction, or $p_1 = \ldots = p_{j-1} = 0$ (if $x_0 \geq i_j$), again a contradiction.\vspace{5mm}

\noindent It follows that there are two constants $0 < c < C < +\infty$ such that $c < \inf_j \inf_r \lambda_j^{(r)} \leq \sup_j \sup_r \lambda_j^{(r)} < C$ and hence it follows that there is a limit point $\underline{\lambda}$ of $\underline{\lambda}^{(r)}$ which provides a solution. The lower bound~\eqref{eqalthcl2} follows by taking the limit as $r\rightarrow +\infty$ in~\eqref{eqalph}. \qed 

\paragraph{Note} At this point there is a (minor) divergence when one tries to establish existence of $m$ such that the generator defined by~\eqref{eqprob2} has the required properties. When considering this problem, there are killing rates $k_1$ and $k_M$ on sites $i_1$ and $i_M$ respectively, which are not necessarily $0$. But after the process reaches either of these sites, the killing rate is exponential and therefore the process survives with positive probability for any finite time and it is straightforward to obtain an upper bound on the killing probability which is strictly less than $1$  when the process is stopped at a Gamma$\left(r, \frac{r}{t}\right)$ time for fixed $t > 0$; the upper bound is independent of $r \geq 1$.

\subsection{Fixed Point Theorem}  For fixed   $\underline{k}$   let $\underline{h} :  \mathbb{R}^M  \times \mathbb{R}^{M-2}_+ \rightarrow \mathbb{R}^M$ denote the function defined by: 

\begin{equation}\label{eqhdef} \underline{h}(\underline{p},  \underline{\lambda}) = \frac{1}{\sum_{j,k} p_j {\cal N}_{jk}(\underline{\lambda})}\underline{p} {\cal N}(\underline{\lambda}) = \frac{1}{\sum_j p_j (1 + \frac{t}{r} \lambda_j \widetilde{k}_j)} \underline{p} {\cal N}(\underline{\lambda}).
\end{equation} 

\noindent where $\widetilde{\underline{k}}$ is defined by~\eqref{eqtildek}. Directly from the definition, for any $\underline{p} \in \mathbb{R}^M_+$ and $\underline{\lambda} \in \mathbb{R}^{M-2}_+$, 

\begin{equation}\label{eqsumh} \sum_j h_j(\underline{p}, \underline{\lambda}) = 1. \end{equation} 

\noindent From the definition, it is also clear that $\underline{h}(\alpha\underline{p}, \underline{\lambda}) = \underline{h}(\underline{p}, \underline{\lambda})$ for any $\underline{p} \in \mathbb{R}^M_+$, $\alpha \in \mathbb{R} \backslash \{0\}$ and $\underline{\lambda} \in \mathbb{R}_+^{M-2}$. \vspace{5mm}

\noindent For $r \geq 2$, set

\begin{equation}\label{eqhrdef} \underline{h}^{(r)}(\underline{p}, \underline{\lambda}) = \underline{h}(\underline{h}^{(r-1)}(\underline{p}, \underline{\lambda}), \underline{\lambda}).
\end{equation}

\begin{Th}[Fixed Point Theorem] \label{thfpt} Set

\begin{equation}\label{eqadef2} {\cal A} (\underline{\lambda}, \underline{p})(j):=  {\cal G}_j \left(\frac{t}{r}, Q \left (\frac{t}{r}, h^{(r-1)}(\underline{p} , \underline{\lambda}), \underline{\lambda}\right )  \right ) \qquad j = 2, \ldots, M-1
\end{equation} 

\noindent  where $Q$ defined by~\eqref{eqqespj2} and ${\cal G}$ by~\eqref{eqgeee}. There exists a solution $\underline{\lambda}$ to the equation

\begin{equation}\label{eqfpeq} \underline{\lambda} = {\cal A}(\underline{\lambda}, \underline{p}).
\end{equation}

\noindent which satisfies $\underline{\lambda} \in \mathbb{R}^{M-2}_+$. 
\end{Th}

\paragraph{Proof of Theorem~\ref{thfpt}} The proof   follows the lines of~\cite{N1}.  As in~\cite{N1}, for $\underline{p} \in \mathbb{R}^M$, set

\begin{equation}\label{eqcee} C(\underline{p}, \epsilon) = \sum_{j=1}^M \left(\frac{p_j}{\sum_{k=1}^M (p_k \vee \epsilon)} \vee \epsilon \right).\end{equation}

\noindent For any $\underline{p} \in \mathbb{R}^M$ and $\epsilon \in (0,1)$, $C(\underline{p}, \epsilon) \leq M$. For $\underline{p} \in \mathbb{R}^M$ such that $\sum_{k=1}^M (p_k \vee 0) \geq 1$, it follows that for any $\epsilon \in [0,1)$, $\sum_{k=1}^M (p_k \vee \epsilon) \geq 1$ and hence that

\[M \geq C(\underline{p},\epsilon) \geq \sum_{j=1}^M \left(\frac{p_j}{\sum_{k=1}^M (p_k \vee \epsilon)} \vee \frac{\epsilon}{\sum_{k=1} (p_k \vee \epsilon) } \right) = 1.\]

\noindent Let ${\cal P}^{(\epsilon)}: \mathbb{R}^M \rightarrow \mathbb{R}_+^M$ denote the function

\begin{equation}\label{eqp} {\cal P}^{(\epsilon)}_j(\underline{p}) = \frac{1}{C(\underline{p}, \epsilon)}\left(\frac{p_j}{\sum_{k=1}^M (p_k  \vee \epsilon)} \vee \epsilon \right) \end{equation}

\noindent where $C$ is defined by~\eqref{eqcee} so that $\sum_{j=1}^M {\cal P}_j^{(\epsilon)}(\underline{p}) = 1$. It follows that for any $\underline{p} \in \mathbb{R}^M$,

\[  \min_j {\cal P}_j^{(\epsilon)}(\underline{p}) \geq \frac{\epsilon}{M}.\]

\noindent Set

\begin{equation}\label{eqaeps2} {\cal A}^{(\epsilon)}(\underline{\lambda}, \underline{p})(j) =  {\cal G}_j\left ( \frac{t}{r},  Q \left (\frac{t}{r}, {\cal P}^{(\epsilon)}(h^{(r-1)}(\underline{p}, \underline{\lambda})), \underline{\lambda} \right )  \right ). \end{equation}

\begin{Lmm} For each $\epsilon > 0$, there exists a $K(\epsilon) < +\infty$  such that 
 \[  \sup_{\underline{\lambda} \in \mathbb{R}^{M-2}_+} \max_{j \in \{2, \ldots, M-1\}} {\cal A}^{(\epsilon)}(\underline{\lambda}, \underline{p})(j) \leq K (\epsilon).\]
\end{Lmm}

\paragraph{Proof} Consider Equation~\eqref{eqgeee}. If $p_j \geq \epsilon$ for all $j \in \{1, \ldots, M\}$, then ${\cal A}^{(\epsilon)}(\underline{\lambda}, \underline{p})(j) \leq f_j$ where $f_j$ satisfies

\[ \left\{ \begin{array}{ll} f_j = a + b\sum_{i=1}^{j-1} f_i & j = 2, \ldots, M \\ f_1 = 0 & \\
a = \frac{1}{\epsilon} \max_j\frac{r(\kappa_{j+1} - \kappa_{j-1})(\kappa_M - \kappa_1)}{t(\kappa_{j+1} - \kappa_j)(\kappa_j - \kappa_{j-1})}, & b = a t \max_j \widetilde{k}_j
\end{array} \right. \]

\noindent where $\widetilde{k}_j : j = 1, \ldots, M$ is defined by~\eqref{eqtildek}. The solution to this equation is

\[ f_1 = 0 \qquad f_j = a(1 + b)^{j-2} \qquad j \geq 2.\]

\noindent and hence

\[ 0 \leq \min_j {\cal G}_j \leq \max_j {\cal G}_j \leq a (1 + b)^{M-2}.\]

\noindent This depends on $\epsilon$, but it does not depend on $\underline{\lambda}$. \qed \vspace{5mm}

\noindent It is clear from the construction that for fixed $\epsilon > 0$, ${\cal A}^{(\epsilon)}(., \underline{p})$ is continuous in $\underline{\lambda}$. Therefore, by the Schauder Fixed Point Theorem, there is a solution $\underline{\lambda}^{(\epsilon)}$ to the equation

\[ \underline{\lambda} = {\cal A}^{(\epsilon)}(\underline{\lambda}, \underline{p}).\]

\noindent Let $\underline{\lambda}^{(\epsilon)}$ denote a fixed point (solution) and let

\begin{equation}\label{eqheps} \underline{h}_\epsilon = {\cal P}^{(\epsilon)}(h^{(r-1)}(\underline{p}, \underline{\lambda}^{(\epsilon)})),
\end{equation}

\noindent where ${\cal P}^{(\epsilon)}$ is defined by~\eqref{eqp}, so that

\begin{equation}\label{eqlamep} \lambda^{(\epsilon)}_j  =  {\cal G}_j \left(\frac{t}{r}, Q \left(\frac{t}{r},\underline{h}_\epsilon \right )  \right ) \qquad j = 2, \ldots, M-1
\end{equation}

\noindent where ${\cal G}$ is defined by~\eqref{eqgeee} and $Q$ is defined by~\eqref{eqqespj2}. It is required to show:

\begin{itemize}
\item $\limsup_{\epsilon \rightarrow 0}\max_{j \in \{2, \ldots, M-1\}} \lambda_j^{(\epsilon)} < +\infty$
\item $\liminf_{\epsilon \rightarrow 0} \min_{j \in \{2, \ldots, M-1\}} \lambda_j^{(\epsilon)} > 0$
\end{itemize}

\paragraph{Showing $\limsup_{\epsilon \rightarrow 0} \max_{j \in \{2, \ldots, M-1\}} \lambda_j^{(\epsilon)} < +\infty$} It follows from~\eqref{eqlamep}, using the definition of ${\cal G}$ (Equation~\eqref{eqgeee}) and the definition of $Q$ (Equation~\eqref{eqqespj2}) that: 

\begin{eqnarray}\lefteqn{ \nonumber\frac{t}{r} h_{\epsilon,j} \lambda^{(\epsilon)}_j \left(1 + \frac{t}{r} \widetilde{k}_j \lambda^{(\epsilon)}_j\right)}\\&& \label{eqhepp2} = \frac{(\kappa_{j+1} - \kappa_{j-1})}{(\kappa_{j+1} - \kappa_j)(\kappa_j - \kappa_{j-1})}\times \left\{ \begin{array}{ll} \sum_{i=1}^{j-1} (\kappa_j - \kappa_i) h_{\epsilon,i} \left( 1 + \frac{t}{r} \widetilde{k}_i \lambda_i^{(\epsilon)}\right)^2 & 2 \leq j \leq l-1 \\ 
\sum_{i=j+1}^{M} (\kappa_i - \kappa_j) h_{\epsilon,i} \left( 1 + \frac{t}{r}\widetilde{k}_i \lambda_i^{(\epsilon)}\right)^2 & l \leq j \leq M -1
\end{array}\right.  
\end{eqnarray}

\noindent It follows that 

\begin{equation}\label{eqhepconv} h_{\epsilon,j} \stackrel{\lambda_j^{(\epsilon)} \rightarrow +\infty}{\longrightarrow} 0.
\end{equation}

\noindent This can be seen inductively from~\eqref{eqhepp2}: recall that $h_{\epsilon,j} > 0$ for each $j$ and $\sum_{j=1}^M h_{\epsilon,j} = 1$. Since $\lambda_1^{(\epsilon)} = \lambda_M^{(\epsilon)} = 0$, the result is clearly true for $j = 2$ and $M-1$ and, furthermore, there are uniform bounds on  $\frac{t}{r} h_{\epsilon,2}\lambda_2^{(\epsilon)}\left(1 + \frac{t}{r}\widetilde{k}_2 \lambda_2^{(\epsilon)} \right )$ and $\frac{t}{r} h_{\epsilon,M-1}\lambda_{M-1}^{(\epsilon)}\left(1 + \frac{t}{r}\widetilde{k}_{M-1} \lambda_{M-1}^{(\epsilon)} \right )$. From this, it follows that there are uniform bounds on $\frac{t}{r} h_{\epsilon,2} \left(1 + \frac{t}{r}\widetilde{k}_2 \lambda_2^{(\epsilon)} \right )^2$ and  $\frac{t}{r} h_{\epsilon,M-1} \left(1 + \frac{t}{r}\widetilde{k}_{M-1} \lambda_{M-1}^{(\epsilon)} \right )^2$. Inductively, it follows  that there are uniform bounds on $\frac{t}{r} h_{\epsilon,j} \left(1 + \frac{t}{r}\widetilde{k}_j \lambda_j^{(\epsilon)} \right )^2$ which hold for all $j$ and hence~\eqref{eqhepconv} follows from~\eqref{eqhepp2}. \vspace{5mm}

\noindent Set 

\[ K(\underline{\lambda}, \epsilon) := \sum_{j=1}^M \left ( h_j^{(r-1)}(\underline{p}, \underline{\lambda}) \vee \epsilon \right)\qquad \mbox{and} \qquad K_\epsilon := K(\underline{\lambda}^{(\epsilon)},\epsilon).\]

\noindent From~\eqref{eqsumh}, it follows that $K_\epsilon \geq 1$. Set

\[ C_\epsilon = C \left (\underline{h}^{(r-1)}(\underline{p}, \underline{\lambda}^{(\epsilon)}) \right )\]
 
\noindent where $C$ is the function defined by~\eqref{eqcee}.\vspace{5mm} 

\noindent Let 

\[ \widetilde{{\cal N}}^{(\epsilon)} = \widetilde{\cal N}(\underline{\lambda}^{(\epsilon)}), \qquad {\cal N}^{(\epsilon)} = {\cal N}(\underline{\lambda}^{(\epsilon)}).\] 

\noindent The first equality below follows from the definition of $\underline{h}_\epsilon$ by~\eqref{eqheps} and~\eqref{eqp}. The second equality follows from the definition of $\underline{h}$ by~\eqref{eqhdef} and $\underline{h}^{(r)}$ given by~\eqref{eqhrdef} together with the identity: $\sum_{j=1}^M h_j^{(r-1)} = 1$, which follows from~\eqref{eqhrdef} and~\eqref{eqsumh}. Recall that $\underline{p}$ and $\underline{h}$ are taken as {\em row} vectors.
 
 \[ h_{\epsilon,j} = \frac{1}{C_\epsilon}\left(\frac{1}{K_\epsilon}h^{(r-1)}_j(\underline{p}, \underline{\lambda}^{(\epsilon)}) \vee \epsilon \right ) = \frac{1}{C_\epsilon K_\epsilon} \frac{1}{\sum_{k} (\underline{p} ({\cal N}^{(\epsilon)})^{r-1})_{k}} (\underline{p} ({\cal N}^{(\epsilon)})^{r-1})_{j} \vee \frac{\epsilon}{C_\epsilon} \]

\noindent where $(\underline{p} ({\cal N}^{(\epsilon)})^{r-1})_{j} = \sum_{k} p_k (({\cal N}^{(\epsilon)})^{r-1})_{k,j}$, $(({\cal N}^{(\epsilon)})^{r-1})_{k,j}$ being the $(k,j)$ component of the matrix $({\cal N}^{(\epsilon)})^{r-1}$. \vspace{5mm}  

\noindent Set 

\begin{equation}\label{eqhhat} \widehat{{\cal H}}_{\epsilon;.,j} = 
\left\{ 
\begin{array}{lll}
\frac{1}{C_\epsilon K_\epsilon \sum_{k} (\underline{p} ({\cal N}^{(\epsilon)})^{r-1})_{k}} (({\cal N}^{(\epsilon)})^{r-1})_{.j} & h_{\epsilon,j} > \frac{\epsilon}{ C_\epsilon}      \\ \frac{\epsilon}{C_\epsilon} & h_{\epsilon,j} = \frac{\epsilon}{ C_\epsilon} & \end{array} \right. 
\end{equation}

\noindent then, since $\sum_{j=1}^M p_j = 1$,

\begin{equation}\label{eqhep} \underline{h}_{\epsilon} = \underline{p} \widehat{{\cal H}}_{\epsilon}.\end{equation}

\noindent Define $\underline{p}^{(\epsilon)}$ as: 

\begin{equation}\label{eqpdef} \underline{p}^{(\epsilon)} := \frac{1}{\sum_{l} ( \underline{p} \widehat{{\cal H}}_{\epsilon} ({\cal N}^{(\epsilon)})^{-(r-1)})_{l}} \underline{p} \widehat{\cal H}_{\epsilon} ({\cal N}^{(\epsilon)})^{-(r-1)}.
\end{equation}

\noindent By construction, $\sum_{j=1}^M p_j^{(\epsilon)} = 1$. Furthermore, it follows from the definition that $\underline{p}^{(\epsilon)}$ satisfies:

\begin{equation}\label{eqpdef2} \underline{p}^{(\epsilon)} = \frac{1}{\sum_{l} ( \underline{h}_{\epsilon} ({\cal N}^{(\epsilon)})^{-(r-1)})_{l}} \underline{h}_{\epsilon}({\cal N}^{(\epsilon)})^{-(r-1)}.
\end{equation}

\noindent From the characterisation given by Lemma~\ref{lmmnprep}, it follows that $0 \leq ((\widetilde{\cal N}^{(\epsilon)})^{-(r-1)})_{j,k} \leq 1$   for each  $(j,k) \in \{1, \ldots, M+1\}^2$. Furthermore, $h_{\epsilon;j} \geq 0$ for all $\epsilon > 0$ and all $j \in \{1, \ldots, M\}$. From this, it follows that $p^{(\epsilon)}_j \geq 0$ for each $j \in \{1, \ldots, M\}$.   \vspace{5mm} 

\noindent From~\eqref{eqpdef} and~\eqref{eqpdef2}, it follows that  

\[ \underline{h}_{\epsilon} = \left( \sum_{l} \left (\underline{p} \widehat{{\cal H}}_{\epsilon} ({\cal N}^{(\epsilon)})^{-(r-1)}\right )_l \right) \underline{p}^{(\epsilon)}({\cal N}^{(\epsilon)})^{r-1}.\]

\noindent  Set 

\begin{equation}\label{eqsdef} S_\epsilon = \left \{ \beta|h_{\epsilon, \beta} = \frac{\epsilon}{C_\epsilon} \right \}
\end{equation}

\noindent Let $Y$ be a continuous time Markov chain with state space $\{1, \ldots, M+1\}$  with transition intensity matrix given by Equation (\ref{eqthetmat}), Definition~\ref{defintmat}. Let $T$ denote an independent time with distribution  $T \sim \mbox{Gamma} \left (r-1,\frac{r}{t}\right )$ (using the notation of~\eqref{eqgamden}). Let 
\[ c(m) = 1 - \mathbb{P}(Y_T = M+1 | Y_0 = m).\]

\noindent It follows from~\eqref{eqtmtnm1rep} that $((\widetilde{{\cal N}}^{(\epsilon)})^{-(r-1)})_{M+1,j} = 0$ for $j = 1, \ldots, M$. From this it follows that 

\[ (\widetilde{{\cal N}}^{(\epsilon)})^{-(r-1)} = \left(\begin{array}{c|c} ({\cal N}^{(\epsilon)})^{-(r-1)} & - ({\cal N}^{(\epsilon)})^{-(r-1)} \underline{v} \\ \hline \underline{0} & 1 \end{array}\right )\]

\noindent where $\underline{0}$ is an $M$-row vector of $0$s, and $\underline{v}$ is the $M$-column vector with $v_j = (\widetilde{{\cal N}}^{(\epsilon) r-1 })_{j,M+1}$, $j = 1, \ldots, M$. 
It follows that for $m_2 \in S_\epsilon$, 

\[ \sum_k (({\cal N}^{(\epsilon)})^{-(r-1)})_{m_1,k}\widehat{{\cal H}}_{k,m_2} = \frac{\epsilon}{C_\epsilon} \sum_{k=1}^M (({\cal N}^{(\epsilon)})^{ -(r-1)})_{m_1,k} = \frac{\epsilon}{C_\epsilon}(1 - ((\widetilde{\cal N}^{(\epsilon)})^{-(r-1)})_{m_1,M+1}) = \frac{\epsilon}{C_\epsilon} c(m_1).\] 

\noindent It follows that: 

\begin{equation}\label{eqnmn} \left(({\cal N}^{(\epsilon)})^{-(r-1)} \widehat{{\cal H}}_\epsilon \right)_{m_1,m_2} = \left\{ \begin{array}{ll}\frac{1}{C_\epsilon K_\epsilon\sum_{jk} p_j (({\cal N}^{(\epsilon)})^{ r-1})_{j,k}} I(m_1,m_2) & m_2 \not \in S_\epsilon \\ \frac{\epsilon}{C_\epsilon}  c(m_1) & m_2 \in S_\epsilon \end{array}\right. \end{equation}

\noindent where $I(m_1,m_2) = \left\{ \begin{array}{ll} 1 & m_1 = m_2 \\ 0 & m_1 \neq m_2  \end{array}\right.$. Set 
 
\begin{equation}\label{eqefff} F_\epsilon := ({\cal N}^{(\epsilon)})^{-(r-1)}  \widehat{{\cal H}}_\epsilon \qquad \mbox{so that} \qquad ({\cal N}^{(\epsilon)})^{r-1} F_\epsilon = \widehat{{\cal H}}_\epsilon. \end{equation}
 
\noindent It follows (from~\eqref{eqpdef}, using~\eqref{eqhep} in the denominator) that: 

\begin{equation}\label{eqpprerewrite} \underline{p}^{(\epsilon)} =  \frac{1}{\sum_{j,k } h_{\epsilon,j} (({\cal N}^{(\epsilon)})^{ -(r-1)})_{j,k}} \underline{p} {\cal N}^{(\epsilon) r-1} F_\epsilon ({\cal N}^{(\epsilon)})^{ -(r-1)}.
\end{equation}

\noindent  Let $\Lambda_\epsilon$ be the matrix such that 

\[ \Lambda_{\epsilon; m_1,m_2} = \left\{ \begin{array}{ll} 1 & m_1 = m_2 \qquad m_2 \not \in S_\epsilon \\ 0 & \mbox{otherwise}. \end{array}\right. \] 

\noindent Let ${\cal I}_\epsilon$ denote the matrix with entries: 
 
\[ {\cal I}_{\epsilon; m_1,m_2} = \left\{\begin{array}{ll} c(m_1) & m_2 \in S_\epsilon \\ 0 & \mbox{otherwise} \end{array}\right. \]
 
\noindent Then ${\cal I}_\epsilon$ has column  $(c(1), \ldots, c(M))^t$  for each $m_2 \in S_\epsilon$ and the remaining columns are columns of $0$s. Then~\eqref{eqnmn} may be written, using  $F_\epsilon$ from~\eqref{eqefff} as:
 
 \[F_\epsilon = \frac{1}{C_\epsilon K_\epsilon \sum_{k} (\underline{p} ({\cal N}^{(\epsilon)})^{r-1})_{k}}\Lambda_\epsilon + \frac{\epsilon}{C_\epsilon }  {\cal I}_\epsilon, \]
 
 \noindent so that

\begin{equation}\label{eqppostrewrite}
 \underline{p}^{(\epsilon)} = \frac{1}{\sum_k \left(\underline{h}_\epsilon ({\cal N}^{(\epsilon)})^{-(r-1)} \right )_k }\left( \frac{\underline{p} ({\cal N}^{(\epsilon)})^{r-1} \Lambda_\epsilon ({\cal N}^{(\epsilon)})^{-(r-1)}}{C_\epsilon K_\epsilon \sum_k (\underline{p} ({\cal N}^{(\epsilon)})^{r-1})_k} + \frac{\epsilon}{C_\epsilon} ({\cal N}^{(\epsilon)})^{r-1} {\cal I}_\epsilon ({\cal N}^{(\epsilon)})^{-(r-1)}\right).
\end{equation}

\noindent Note that, since $((\widetilde{{\cal N}}^{(\epsilon)})^{ -(r-1)})_{jk} = (({\cal N}^{(\epsilon)})^{-(r-1)})_{jk}$ for $(j,k) \in \{1, \ldots, M\}^2$, and \[ \sum_{k=1}^{M+1} ((\widetilde{{\cal N}}^{(\epsilon)})^{-(r-1)})_{jk} = 1 \qquad \forall j \in \{1, \ldots, M\},\] it follows that

\begin{equation}\label{eqhnk} \sum_k \left(\underline{h}_\epsilon ({\cal N}^{(\epsilon)})^{-(r-1)} \right )_k = 1 - \sum_k h_{\epsilon;k} ((\widetilde{{\cal N}}^{(\epsilon)})^{-(r-1)})_{k,M+1} > c_1 > 0
\end{equation}

\noindent for a $c_1 > 0$ which does not depend on $\epsilon$, by Lemma~\ref{lmmcruclb}.

\paragraph{Note} Similarly to the note at the end of Subsection~\ref{mod}, this is the other point where an additional argument is required when an infinitesimal generator given by~\eqref{eqprob2}, since there is killing at sites $i_1$ and $i_M$ at rates $k_1$ and $k_M$ respectively, which are not necessarily $0$. The modification is similar. Consider the proof of Lemma~\ref{lmmcruclb}; after the process eventually reaches state $1$ or state $M$, which it does with positive probability, the killing rate after it hits these sites is bounded; it has rate $k_1$ on site $1$ and $k_M$ on site $M$ and with Generator~\eqref{eqprob2} the kill rate does not depend on $\underline{\lambda}$. Hence a $c_1 > 0$ may be obtained independent of $\epsilon$ such that~\eqref{eqhnk} holds. \vspace{5mm}

\noindent For any invertible matrix ${\cal S}$, the eigenvalues of ${\cal S}^{-1}{\cal A} {\cal S}$ are the same as the eigenvalues of ${\cal A}$. It follows that the eigenvalues of $({\cal N}^{(\epsilon)})^{r-1} \Lambda_\epsilon ({\cal N}^{(\epsilon)})^{-(r-1)}$ are the eigenvalues of $\Lambda_\epsilon$; $0$ with multiplicity equal to the number of elements of $S_\epsilon$ and the remaining eigenvalues all $1$. Similarly, the eigenvalues of ${\cal I}_\epsilon$ are bounded independently of $\epsilon$, since each entry of the $M\times M$ matrix lies in $[0,1]$. It follows that 

\[ \lim_{\epsilon \rightarrow 0}\frac{\epsilon}{C_\epsilon} \underline{p} ({\cal N}^{(\epsilon)})^{r-1}  {\cal I}_\epsilon ({\cal N}^{(\epsilon)})^{ -(r-1)} = 0.\]   \vspace{5mm}

\noindent It now follows directly that if $K_\epsilon  \sum_{k} (\underline{p} ({\cal N}^{(\epsilon)})^{r-1})_{k}   \stackrel{\epsilon \rightarrow 0}{\longrightarrow} + \infty$, then $\underline{p}^{(\epsilon)}\stackrel{\epsilon \rightarrow 0}{\longrightarrow} \underline{0}$, contradicting the fact that $p_j^{(\epsilon)} \geq 0$ for each $j$ and $\sum_j p_j^{(\epsilon)} = 1$ for each $\epsilon \in (0,1)$. \vspace{5mm}
 
\noindent Therefore: 

\[ 0 \leq \inf_\epsilon K_\epsilon \left(\sum_{k} (\underline{p} ({\cal N}^{(\epsilon)})^{r-1})_{k} \right) \leq \sup_\epsilon K_\epsilon \left(\sum_{k} (\underline{p} ({\cal N}^{(\epsilon)})^{r-1})_{k} \right) < +\infty.\]
 
\noindent From the definition of $K_\epsilon$, 

\[ K_\epsilon \left(\sum_{k} (\underline{p}  ({\cal N}^{(\epsilon)})^{ r-1})_{k} \right)  = \sum_{k=1}^M \left( (\underline{p} ({\cal N}^{(\epsilon)})^{r-1})_{k} \vee \left(\sum_{k} (\underline{p} ({\cal N}^{(\epsilon)})^{r-1})_{k}\right)\epsilon \right). \]

\noindent From the above, 
\[ 0 \leq \inf_\epsilon \sum_{k} (\underline{p} ({\cal N}^{(\epsilon)})^{ r-1})_{k} \leq \sup_\epsilon \sum_{k} (\underline{p} ({\cal N}^{(\epsilon)})^{r-1})_{k} < +\infty\]

\noindent and

\[ \sup_\epsilon \max_k \left((\underline{p} ({\cal N}^{(\epsilon)})^{ r-1})_{k} \vee 0 \right) < + \infty,\]

\noindent from which 

\[ \sup_\epsilon \max_k \left | (\underline{p} ({\cal N}^{(\epsilon)})^{ r-1})_{k}   \right |  < + \infty.\]

\noindent Set $\lambda^{*(\epsilon)} = \max_j \lambda_j^{(\epsilon)}$ and let 

\[{\cal N}^{*(\epsilon)} = \frac{1}{\lambda^{*(\epsilon)}}{\cal N}^{(\epsilon)} \]

\noindent (that is, divide every element by $\lambda^{*(\epsilon)}$).  Then if there is a sequence $\epsilon_n \rightarrow 0$ such that $\lambda^{*(\epsilon_n)} \stackrel{n \rightarrow +\infty}{\longrightarrow} + \infty$,  any limit point ${\cal N}^*$ of ${\cal N}^{*(\epsilon_n)}$ satisfies

\[ 0 = \underline{p}{\cal N}^{*(r-1)}.\]

\noindent   It follows from the construction of ${\cal N}^*$ that the rank $\rho$ of ${\cal N}^*$ is the number of components of $\underline{\lambda}$ such that $\lim_{n \rightarrow +\infty}\frac{\lambda_j^{(\epsilon_n)}}{\lambda^{*(\epsilon_n)}} > 0$, where  where $\underline{\lambda}^{(\epsilon_n)}$ is a sequence that gives the limit point. This is seen as follows: consider the lowest index $k_1$ such that $\lim_{n \rightarrow +\infty} \frac{\lambda_{k_1}^{(\epsilon_n)}}{\lambda^{*(\epsilon_n)}} > 0$, then ${\cal N}^*$ 
 in the limit, column $k_1 - 1$ will have exactly one entry; element ${\cal N}^*_{k_1, k_1 - 1}$ will be the only non-zero element of column $k_1$. 
Suppose $k_1 < \ldots < k_\rho$ are the relevant indices, then the columns $({\cal N}^*_{.,k_1-1}, \ldots, {\cal N}^*_{.,k_\rho-1})$ provide an
upper triangular matrix, with elements ${\cal N}^*_{k_j, k_j-1} \neq 0$ and ${\cal N}_{p,k_j - 1} = 0$ for all $p \geq k_j + 1$, proving that
 ${\cal N}^*$ is of rank $\rho$.\vspace{5mm}
 
\noindent Therefore ${\cal N}^{*(r-1)}$ is of rank $\rho$ and the non-zero rows of ${\cal N}^{*(r-1)}$ are those corresponding to
the indices $k : \lim_{n \rightarrow +\infty} \frac{\lambda_k^{(\epsilon_n)}}{\lambda^{*(\epsilon_n)}} > 0$. Since the space spanned by the $\rho$ rows is of rank $\rho$, it follows that
 $p_k = 0$ for each of these $p_k$, which is a contradiction (since, by hypothesis, $p_k > 0$ for each $k$).
Hence

\[ \sup_\epsilon \lambda^{*(\epsilon)} < +\infty.\]

\paragraph{Showing that $\inf_\epsilon \min_{j \in \{2, \ldots, M-1\}} \lambda_j^{(\epsilon)} > 0$.}  Now suppose that there is a subsequence $\lambda_j^{(\epsilon_n)} \stackrel{n \rightarrow +\infty}{\longrightarrow} 0$ for some $j \in \{2, \ldots, M-1\}$. As before, $\lambda^{*(\epsilon)} = \max_{j \in \{2, \ldots, M-1\}}  \lambda^{(\epsilon)}_j$. Recall the representation from  Lemma~\ref{lmmnprep}, that 

\[ (\widetilde{{\cal N}}^{-(r-1)})_{j,k} (\underline{\lambda}) = \mathbb{P}(Y_T = k | Y_0 = j)\]

\noindent where $Y$ is a continuous time Markov chain with state space $\{1, \ldots, M+1\}$, with intensity matrix given by Equation~\eqref{eqthetmat}, Definition~\ref{defintmat} and $T$ is an independent random variable with distribution 
$T \sim \mbox{Gamma}(r-1,\frac{r}{t})$ (using parametrisation found in~\eqref{eqgamden}).  Let $\pi_T$ denote the density function of the random variable $T$. Recall that $\sup_\epsilon \max_k \lambda^{(\epsilon)}_k < +\infty$ and suppose that  $\lambda_j^{(\epsilon_n)} \rightarrow 0$

If $\lambda^{(\epsilon_n)}_j \rightarrow 0$ for some $m_1 \leq j < m_2$ where $m_1 < m_2$, then, letting $\tau_j = \inf\{r | X_r = j\}$ and $\pi_j (dr) $ the probability measure such that $\mathbb{P}(\tau_j \in A) = \int_A \pi_j (dr)$, then  

\begin{eqnarray*}\lefteqn{ \mathbb{P}(Y_T^{(\epsilon_n)} =  m_2 |Y_0^{(\epsilon_n)} =  m_1 ) = \int_0^\infty \mathbb{P}(Y_s^{(\epsilon_n)} =  m_2|Y_0^{(\epsilon_n)} = m_1, T = s)\pi_T(s) ds}\\&&
 \int_0^\infty \int_0^s \mathbb{P}(Y_{r}^{(\epsilon_n)} = j |Y_0^{(\epsilon_n)} =  m_1)\mathbb{P}(Y_{s - r}^{( \epsilon_n)} = m_2 | Y_0^{(\epsilon_n)} = j) \pi_j(dr) \pi_T(s)ds 
\end{eqnarray*}

\noindent so that if $\lambda_{j}^{(\epsilon_n)} \rightarrow 0$, then $\mathbb{P}(Y_s^{(\epsilon_n)} = k|Y_0^{(\epsilon_n)} = j) \rightarrow 0$ for all $k$. It follows that

\begin{equation}\label{eqytepn} \mathbb{P}(Y_T^{(\epsilon_n)} =  m_2 |Y_0^{(\epsilon_n)} = m_1) = ({\cal N}^{-(r-1)})_{m_1,m_2}(\underline{\lambda}^{(\epsilon_n)}) \rightarrow 0
\end{equation}

\noindent for all $(m_1, m_2)$ such that $m_1 \leq j < m_2$. Similarly, if $\lambda^{(\epsilon_n)}_j \stackrel{n \rightarrow +\infty}{\longrightarrow} 0$  for some $m_1 \geq j > m_2$ where $m_1 > m_2$, then~\eqref{eqytepn} holds.\vspace{5mm}

\noindent It follows that $ {\cal N}^{-(r-1)}_{kp} (\underline{\lambda}^{(\epsilon_n)})  \rightarrow 0$ for
all $(k,p)$ such that $k \leq j < p$ or $k \geq j > p$.   \vspace{5mm}

\noindent Furthermore, it follows from~\eqref{eqhepp2}  that for any sequence with limit point $\underline{\lambda}^{(0)}$ such that $\lambda^{(\epsilon_n)}_j\rightarrow 0$ for some $j \in \{2, \ldots, M-1\}$, there is an $l \in \{1, \ldots, M\}$ such that if $j \leq l-1$, then $h_k^{(r-1)}(\underline{p}, \underline{\lambda}^{(0)}) \leq 0$ for all $1 \leq k \leq j - 1$ and if $j \geq l$ then $h_k^{(r-1)} (\underline{p}, \underline{\lambda}^{(0)}) \leq 0$ for all $j+1 \leq k \leq M$.\vspace{5mm}

\noindent Now recall the definition of $h^{(r-1)}$ (~\eqref{eqhdef} and~\eqref{eqhrdef}) from which it follows that

\[ h^{(r-1)}(\underline{p},\underline{\lambda}) = \left(\sum_{jk} p_j \left({\cal N}^{r-1}(\underline{\lambda})\right)_{jk} \right)^{-1} \underline{p} {\cal N}^{r-1}(\underline{\lambda}).\]

\noindent From this it follows that: 

\begin{equation} \underline{p} = \underline{p}{\cal N}^{(r-1)}(\underline{\lambda}) {\cal N}^{-(r-1)} (\underline{\lambda}) = \left ( \sum_{j,k} p_j ({\cal N}^{r-1})_{j,k}(\underline{\lambda})\right) \underline{h}^{(r-1)}(\underline{p}, \underline{\lambda}){\cal N}^{-(r-1)}(\underline{\lambda}).
\end{equation}

\noindent Recall the definition of $\lambda^{(\epsilon)}$ given by~\eqref{eqlamep} and let $l$ denote the index from the definition of ${\cal G}$ in~\eqref{eqgeee}. With $\underline{\lambda} = \underline{\lambda}^{(0)}$ and considering the zeroes of ${\cal N}^{-(r-1)}(\underline{\lambda}^{(0)})$, it follows that if $j \leq l-1$, then $p_1 \leq 0, \ldots, p_{j-1} \leq 0$, which is a contradiction. If $j \geq l$, then $p_{j+1} \leq 0, \ldots, p_M \leq 0$, which is a contradiction.\vspace{5mm}

\noindent It follows that any limit point $\underline{\lambda}$ satisfies $0 < \min_{j \in \{2, \ldots, M-1\}} \lambda_j \leq \max_{j \in \{2, \ldots, M-1\}} \lambda_j < +\infty$,  consequently that $h_{0,1} > 0, \ldots, h_{0,M} > 0$, therefore $\underline{h}(\underline{p}, \underline{\lambda}) = \underline{h}_0$ and therefore $\underline{\lambda}$ satisfies Equation~\eqref{eqfpeq}. Theorem~\ref{thfpt}  is proved. \qed

\section{Proof of Theorem~\ref{thdrk}}

 Following the proof of Theorem~\ref{thctlim2}, the theorem is already proved for a finite state space ${\cal S} = \{i_1, \ldots, i_M\}$; let $\underline{a}$ satisfy $a_j =  \lambda_j (i_{j+1} - i_j)(i_j - i_{j-1})$ then, following Lemma~\ref{lmmfininfgen2}, the continuous time, time homogeneous Markov process $X$ that satisfies Theorem~\ref{thctlim2} has infinitesimal generator

\[ {\cal L} = a\left(\frac{1}{2}\Delta + b\nabla - k\right )  \]

\noindent where the operators $\Delta$ and $\nabla$ are defined by Equations~\eqref{laplace} and~\eqref{nabla} respectively, where ${\cal L}$ means:

\[ {\cal L}f(i_j) = a_j \left(\frac{1}{2}\Delta + b_j \nabla - k_j \right) f(i_j) \qquad j = 1, \ldots, M \qquad a_1 = a_M = 0.\]

\noindent For a probability distribution $\mu$ over $\mathbb{R}$, the proof follows the same lines as the proof already given for $k \equiv 0$.  Set

\begin{equation}\label{eqesen} {\cal S}_N = \{i_{N,1}, \ldots, i_{N,M_N}\} \end{equation}

\noindent and let $\underline{p}^{(N)}$ be defined by Equation~\eqref{eqpenj}.  Let $\underline{\lambda}_N$ denote a solution to the terminal distribution problem for distribution $\underline{p}^{(N)}$ over space ${\cal S}_N$. Let 

\[a^{(N)}_j =  a^{(N)}(i_{N,j}) =  \lambda^{(N)}_j (i_{N,j+1} - i_{N,j})(i_{N,j} - i_{N,j-1}) \qquad j = 2, \ldots, M_N - 1  \]

\noindent and let ${\cal L}^{(N)}$ be the infinitesimal generator defined by

\[ {\cal L}^{(N)} f(i_{N,j}) = a^{(N)}(i_{N,j}) \left(\frac{1}{2}\Delta_N + b^{(N)}_j  \nabla_N - k^{(N)}_j \right) f(i_{N,j}) \qquad j = 2, \ldots, M_N - 1\]

\noindent where $\Delta_N$ and $\nabla_N$ are the Laplacian and gradient operators defined on ${\cal S}_N$ (Definition~\ref{deflapder}), the approximate drift field $(b^{(N)}_2, \ldots, b^{(N)}_{M_N-1})$ defined by~\eqref{eqdd} and 

\begin{equation}\label{eqkconv} k^{(N)}_j =  
\frac{1}{i_{N,j+1} - i_{N,j}}\int_{i_{N,j}-}^{i_{N,j+1}} \widehat{k} (x) dx \qquad  j = 2, \ldots, M_N-1  \end{equation}

\noindent where $\widehat{k}$ is from~\eqref{eqtildk}, $\int_{a-}^b$ means integration over the interval $[a,b)$. Then ${\cal L}^{(N)}$ is the infinitesimal generator of the process $X^{(N)}$  with state space ${\cal S}_N \cup \{D\}$, where $D$ denotes a cemetery, such that there is an $l_N$, an $\alpha_N \in (0,1)$ and a $\beta_N \in (0,1)$ such that 

\begin{equation}\label{eqtildyn} \left\{\begin{array}{l} \beta_N \mathbb{P} \left ( X_t^{(N)}  = i_{N,j} | X_0 = i_{N,l_N} \right ) + (1-\beta_N )\mathbb{P} \left ( X_t^{(N)} = i_{N,j} | X_0 = i_{N,l_N -1} \right ) = \alpha_N p^{(N)}_j \\  j = 1, \ldots, M_N \\ \beta_N \mathbb{P} \left ( X_t^{(N)}  \in \{ D \} | X_0 = i_{N,l_N } \right ) + (1-\beta_N )\mathbb{P} \left ( X_t^{(N)} \in \{ D \} | X_0 = i_{N,l_N -1} \right ) = 1 - \alpha_N \end{array}\right.  
\end{equation}

\noindent The quantity $\beta_N$ may be interpreted in  the following way: there is a point $x_{N,0} \in (i_{N,l_N-1}, i_{N,l_N}]$, denoting the initial condition, such that 

\[ \mathbb{P} \left (X_{0+}^{(N)} = i_{N,l_N} | X_0^{(N)} = x_{N,0} \right ) = \beta_N, \qquad \mathbb{P} \left (X_{0+}^{(N)} = i_{N,l_N-1} | X_0^{(N)} = x_{N,0} \right ) = 1 - \beta_N.\]

\noindent Let $\widetilde{X}^{(N)}$ denote the process with infinitesimal generator $a^{(N)}\left(\frac{1}{2}\Delta_N + b^{(N)} \nabla_N\right)$, then
\[ X^{(N)}_t = \left\{\begin{array}{ll} \widetilde{X}^{(N)}_t & t \leq \tau_N \\ D & t > \tau_N \end{array}\right. \]

\noindent where $\tau_N$ is a random time satisfying
\begin{equation}\label{eqtauen} \mathbb{P}(\tau_N \geq s |( \widetilde{X}^{(N)}_.)) = \exp\left\{ -\int_0^s a^{(N)}(\widetilde{X}_r^{(N)})k( \widetilde{X}^{(N)}_r) dr\right \} \qquad s \geq 0.\end{equation}

\noindent It follows from Theorem~\ref{thchgd} that for a sequence $z_j \rightarrow z$, there exists a $\widetilde{X}$ such that 

\begin{equation}\label{eqxlim}\lim_{j \rightarrow +\infty} \mathbb{P} \left( \sup_{0 \leq s \leq t} \left | \widetilde{X}_s^{(N_j)}(z_{N_j}) - \widetilde{X}_s(z) \right | > \epsilon \right ) = 0.\end{equation}

\noindent Let  $\tau$ denote a random time satisfying

\begin{equation}\label{eqrttau} \mathbb{P}(\tau  \geq s |(\widetilde{X}_.)) = \exp\left\{ -\int_0^s \frac{dK}{dm}  (\widetilde{X}_r) dr\right \} \qquad   s \geq 0. \end{equation}

\noindent It follows from Equations (\ref{eqtauen}) and (\ref{eqrttau}) that for a sequence $z_{N_j} \rightarrow z$ such that $z_{N_j} \in {\cal S}_{N_j}$ for each $j$ (${\cal S}_N$ defined by~\eqref{eqesen}), 

\[ \lim_{j \rightarrow +\infty} \left | \mathbb{P}\left ( \left \{ \tau^{(N_j)} \leq t \right \} \right ) - \mathbb{P} \left ( \left \{ \tau \leq t \right \} \right )\right | = 0\]

\noindent and 

\[ \lim_{j \rightarrow +\infty} \sup_{x \in \mathbb{R}} \left | \mathbb{P}(\{\widetilde{X}^{(N_j)}_{t}(z_{N_j}) \leq x\} \cap \{\tau^{(N_j)} > t \}) - \mathbb{P}(\{\widetilde{X}_t(z) \leq x \} \cap \{ \tau > t\})\right | = 0\]
 
\noindent from which it follows that $X$ is a process with infinitesimal generator $\frac{1}{2} \frac{d^2}{dm dx} + \frac{\partial B}{\partial m}\nabla_m - \frac{dK}{dm}$ with the required distribution at the prescribed time $t > 0$, provided $\inf_N \alpha_N > 0$. \vspace{5mm}

\noindent Finally, it has to be shown that for the sequence of measures $m^{(N)}$ there does not exist a subsequence such that $m^{(N_j)}  \rightarrow 0$, which would correspond to $\mathbb{P}(X_t^{(N_j)} \in \{D\}) \stackrel{j \rightarrow +\infty}{\longrightarrow} 1$. \vspace{5mm}

\noindent Let $m^{(N)}$ denote the sequence of measures corresponding to the atomised state spaces. Let $e_{N-}$ and $e_{N+}$ be the indices defined by~\eqref{eqemminemplusN}. Let 
$K_N = m^{(N)}([i_{e_{N-}}, i_{e_{N+}}])$ and $\widehat{m}^{(N)} = \frac{m^{(N)}}{K_N}$. If $K_{N_j} \rightarrow 0$, then there exists a limit point $\widehat{m}$ of $\widehat{m}^{(N_j)}$ such that ${\cal L} := \frac{1}{2}\frac{\partial^2}{\partial \widehat{m} \partial x} + \frac{\partial B}{\partial \widehat{m}} \nabla_{\widehat{m}}   - \frac{\partial K}{\partial \widehat{m}}$ is the infinitesimal generator of a process $X$  which, conditioned on being alive, has stationary distribution $\mu$. 

Let $Y$ denote the process with infinitesimal generator $\frac{1}{2}\frac{\partial^2}{\partial \widehat{m} \partial x} + \frac{\partial B}{\partial \widehat{m}} \nabla_{\widehat{m}}$ and let $P(t;x,dy)$ denote its transition kernel. Let $p(t;x,y) = \frac{P(t;x,dy)}{\widehat{m}(dy)}$. Let $\phi(y) = \frac{\mu(dy)}{\widehat{m}(dy)}$. Then $\phi(y)$ has representation

\begin{equation}\label{eqphirep} \phi(y) = \int \mu (dx) \int p(t;x,y) \frac{\mathbb{E}_x \left [ e^{-\int_0^t \frac{\partial K}{\partial \mu}(Y_s)\phi(Y_s) ds} | Y_t = y \right ]}{\mathbb{E}_x \left [ e^{-\int_0^t \frac{\partial K}{\partial \mu}(Y_s)\phi (Y_s) ds} \right ]}.
\end{equation}

\noindent This can be seen as follows: the transition kernel $Q(t;x,dy)$ for the process $X$ satisfies

\[ \left\{ \begin{array}{l} \frac{\partial}{\partial t} Q(t;x,A) = \frac{1}{2}\frac{\partial^2}{\partial \widehat{m} \partial x} Q(t;x,A) + \frac{\partial B}{\partial \widehat{m}} \nabla_{\widehat{m}} Q(t;x,A) - \frac{\partial K}{\partial \widehat{m}} Q(t;x,A) \\ Q(t;x,A) = {\bf 1}_A(x). \end{array}\right. \]

\noindent Using $\frac{\partial K}{\partial \widehat{m}} = \frac{\partial K}{\partial \mu}\frac{\partial \mu}{\partial \widehat{m}} = \frac{\partial K}{\partial \mu} \phi$, this has representation:

\[ Q(t;x,A) = \mathbb{E}_x \left [ {\bf 1}_A(Y_t)e^{-\int_0^t \frac{\partial K}{\partial \mu}(Y_s)\phi(Y_s) ds}\right ] \qquad A \in {\cal B}(\mathbb{R}).\]

\noindent Conditioning on being alive, $\mathbb{P}(X_t \in A | X_0 = x, X_t \not \in \{D\}) = \frac{Q(t;x,A)}{Q(t;x, \mathbb{R})}$, so that 

\[ \mu (A) = \int \mathbb{P}(X_t \in A | X_0 = x, X_t \not \in \{D\}) \mu (dx) \]

\noindent from which Equation~\eqref{eqphirep} follows. This holds for all $t > 0$. Firstly, it follows from this that $\phi(y) \stackrel{y \rightarrow \pm \infty}{\longrightarrow} 0$. Secondly, by the hypothesis on $b$ and $\mu$, it follows from Lemma~\ref{lmmkappdef} that $Y$ may be put into `martingale' coordinates (described in Section~\ref{subcoc}). Let $z_+ = \sup\{x | x \in \mbox{suppt}(\mu)\}$ and $z_- = \inf\{x | x \in \mbox{suppt}(\mu)\}$.  It follows from basic properties of martingales that if either $z_- > -\infty$ or $z_+ < +\infty$, then $Y_t$ has a well defined limit almost surely, otherwise  $|Y_s| \stackrel{s \rightarrow +\infty}{\longrightarrow} +\infty$. In all cases, it follows from Hypothesis~\ref{hybmu} Part~\ref{hypart4} that $\frac{\partial K}{\partial \mu} (Y_s) \stackrel{s \rightarrow +\infty}{\longrightarrow} 0$. From this, it follows that for fixed $(x,y) \in (z_-,z_+)$, 

\[ \limsup_{t \rightarrow +\infty} \frac{\mathbb{E}_x \left [ e^{-\int_0^t \frac{\partial K}{\partial \mu}(Y_s)\phi(Y_s) ds} | Y_t = y \right ]}{\mathbb{E}_x \left [ e^{-\int_0^t \frac{\partial K}{\partial \mu}(Y_s)\phi (Y_s) ds} \right ]} < 1.\]

\noindent It follows from the existence of a transformation to martingale coordinates that $p(t;x,y) \rightarrow 0$ for all $(x,y) \in (z_-,z_+)$, from which it follows that $\phi(y) \equiv 0$, hence $\mu \equiv 0$ and a contradiction has been obtained.
\qed 

\section{Conclusion and Further Study} The article~\cite{N1} established existence of generalised diffusion to meet a given marginal for any probability measure over $\mathbb{R}$. This article deals with the introduction of drift and killing and establishes conditions on given drift and killing under which there exists a `clock' such that the process, conditioned on being alive at a fixed time $t$, has the prescribed marginal. 

The open problem of interest is to determine the extent to which the conditions on the drift and killing are merely technicalities to make the proofs work, or whether counter examples can be obtained. In particular, can one find a solution to the problem if there exists a string $m$ such that $\frac{1}{2}\frac{\partial^2}{\partial m \partial x} + \frac{\partial B}{\partial m}\nabla_m - \frac{\partial K}{\partial m}$ is the generator of a process which, conditioned on being alive, has {\em invariant measure} $\mu$? This situation (of course) does not arise in the absence of drift and killing. 

Another problem of great interest is to explore the connections between the method given here and the Local Variance Gamma Model by Peter Carr, discussed in~\cite{Carr}. This model considers a process composed with a Gamma process. This boils down to a generalised diffusion stopped at an exponential time. There are further developments in~\cite{Carr} and it is of interest to explore the connections between the process stopped at a Gamma time described here, with the problem of introducing more uniform maturity spacings for the problem of calibrating to meet multiple smiles discussed in~\cite{Carr}.

\setlength{\baselineskip}{2ex}

\end{document}